\documentclass[11pt]{article}

\usepackage{natbib} % for references

\usepackage{color} % coloring ref numbers
\usepackage{hyperref} % coloring url links
\hypersetup{
    colorlinks=true,
    linkcolor=red,
    citecolor=blue,
    urlcolor=blue,
    pdfborder={0 0 0}
}
\usepackage{authblk}        % add authors cleanly
\usepackage[utf8]{inputenc} % allow utf-8 input
\usepackage[T1]{fontenc}    % use 8-bit T1 fonts
\usepackage{url}            % simple URL typesetting
\usepackage{booktabs}       % professional-quality tables
\usepackage{amsfonts}       % blackboard math symbols
\usepackage{microtype}      % micro-typography
\usepackage{lipsum}         % generates dummy text
\usepackage{geometry}       % gives page dimensions
\usepackage{amsmath, amssymb, amsthm} %AMS packages for math stuff
\usepackage{rsfso}      % less slanted version for Ralph Smith Formal Script
\usepackage[cal=cm,        % callligraphy alphabet
scr=rsfso,                   %script alphabet 
calscaled=.96                %alphabet size 
]{mathalfa}                  %math font styles
\usepackage{dsfont}
\usepackage{comment}     % for commenting out sections.

\usepackage{caption}  %used to modify captions for figs, tables, etc.
\usepackage{tikz}
\usepackage{pgfplots}
\pgfplotsset{compat = newest}

\usepackage[toc,page]{appendix}  % For appendices

\newtheorem{proposition} {Proposition} [section]
\theoremstyle{definition}
\newtheorem{definition}{Definition}[section]
\newtheorem{example}{Example}[section]
\theoremstyle{remark}
\newtheorem{remark}{Remark}[section]

\newcommand{\sw}[2][\alpha]{\|{#2}\|_{\psi_{#1}}}

\newcommand{\argmin}{\operatornamewithlimits{argmin}}

\newcommand{\cov}{\mbox{Cov}}

\newcommand{\bs}[1]{\boldsymbol{#1}}

\newcommand{\fr}[2][n]{\frac{{#2}}{{#1}}}
\newcommand{\parenf}[1]{\left(#1\right)}
\newcommand{\parenc}[1]{\left\{#1\right\}}
\newcommand{\parens}[1]{\left[#1\right]}

\newcommand{\frsum}{\fr{1}\sum_{i=1}^n}
\newcommand{\tr}{\text{tr}}

\providecommand{\keywords}[1]{\textbf{\textit{Keywords:}} #1}

\title{Optimal Sparse Estimation of High Dimensional Heavy-tailed Time Series }

\author[1]{Sagnik Halder\thanks{Email:shalder@ufl.edu}}
\author[2]{George Michailidis\thanks{Email: gmichail@ufl.edu}}
\affil[1]{Department of Statistics, University of Florida}
\affil[2]{Department of Statistics \& Informatics Institute, University of hFlorida}

\usepackage{xparse}

\ExplSyntaxOn

\NewDocumentCommand{\definealphabet}{mmmm}
 {% #1 = prefix, #2 = command, #3 = start, #4 = end
  \int_step_inline:nnn { `#3 } { `#4 }
   {
    \cs_new_protected:cpx { #1 \char_generate:nn { ##1 }{ 11 } }
     {
      \exp_not:N #2 { \char_generate:nn { ##1 } { 11 } }
     }
   }
 }

\ExplSyntaxOff

\definealphabet{mb}{\mathbb}{A}{Z}
\definealphabet{mc}{\mathcal}{A}{Z}
\definealphabet{mr}{\mathrm}{A}{Z} % roman/text mode in math env.
\definealphabet{mf}{\mathfrak}{A}{Z}
\definealphabet{mf}{\mathfrak}{a}{z}
\definealphabet{ms}{\mathscr}{A}{Z} %from pkg 'mathrsfs' 
\definealphabet{msf}{\mathsf}{A}{Z} %from pkg 'mathrsfs' 
\definealphabet{msf}{\mathsf}{a}{z}

\begin{document}

\maketitle

\begin{abstract}
Recently, high dimensional vector auto-regressive models (VAR), have attracted   a lot of interest, due to novel applications in the health, engineering and social sciences. The presence of temporal dependence poses additional challenges to the theory of penalized estimation techniques widely used in the analysis of their \textit{iid} counterparts. However, recent work (e.g., \citep{basu2015regularized,kock2015oracle}) has established \textit{optimal} consistency of $\ell_1$-LASSO regularized estimates applied to models involving high dimensional $\textit{stable}$ $\textit{Gaussian}$ processes. The only price paid for temporal dependence is an extra  multiplicative factor that equals 1 for independent and identically distributed (iid) data. Further, \citep{wong2020lasso} extended these results to heavy tailed VAR\textit{s} that exhibit "$\beta$-mixing" dependence, but the rates rates are \textit{sub}-optimal, while the extra factor is intractable.

This paper improves these results in two important directions: (i) We establish \textit{optimal} consistency rates \textit{and} corresponding finite sample bounds for the underlying model parameters that match those for iid data, modulo a price for temporal dependence, that is easy to interpret and equals 1 for iid data. (ii) We incorporate more general penalties in estimation (which are $\textit{not}$ decomposable unlike the $\ell_1$ norm) to induce general sparsity patterns. The key technical tool employed is a novel, easy-to-use $\textit{concentration bound}$ for heavy tailed linear processes, that do $\textit{not}$ rely on "mixing" notions and give tighter bounds.
\end{abstract}


\section{Introduction}

Multivariate time series data are ubiquitous in many application domains, including economics \citep{stock2016dynamic}, finance \citep{geraci2018measuring}, functional genomics \citep{michailidis2013autoregressive} and neuroscience \citep{seth2015granger}. However, statistical models developed for such data often require regularization of their parameters to enable their consistent estimation due to their high dimensional and limited availability of samples. 
%Popular regularized include LASSO and the  "Dantzig" selector  \citep{wainwright2019high}.

We start our exposition by focusing on the stochastic regression given by
\begin{gather}
y_t=x_t^\top\beta^*+\epsilon_t,\quad t=1,\dots,n,    
\end{gather}
where the responses $\{y_t\}$ depend on strictly stationary, centered stochastic processes $\{x_t\}$ in $\mbR^p$ (the predictors), and noise $\{\epsilon_t\}$ in $\mbR$. %\textcolor{red}{the stochastics are the X and the noise, not the y's}
The goal is to consistently estimate the $p$-dimensional regression parameter $\beta^*$ from observed data $\{(y_t,x_t),1\le t\le n\}$, under a high-dimensional regime $p\gg n$. It is common to assume that $\beta^*$ has an underlying low-dimensional structure, such as being sparse or group sparse \citep{wainwright2019high}, which is induced using a penalized estimation approach. 
%Then we derive non-asymptotic consistency and \textit{finite} sample bounds and apply these results to high dimensional autoregressive models driven by SubWeibull noise. 

Two popular penalized methods have been extensively studied in this regard --- the $\textit{Dantzig}$ estimator \eqref{dantzig-selector-def} and $\textit{LASSO}$ \eqref{lasso-def}. There has also been work on more general penalty norms that have proved useful in a number of applications.  For example,  \citet{bach2011optimization,van2014weakly, negahban2012unified} consider norms that are $\textit{weakly decomposable}$, while  \citet{banerjee2015estimation} consider $\textit{atomic}$ norms (defined in section \eqref{sec-model-and-background}).

A large body of literature regarding the statistical properties of regularized estimators for $\beta^*$ exists (see, e.g., \citep{wainwright2019high}), assuming either a fixed design for the predictors or that they are iid samples from some (sub)-Gaussian distribution. In the presence of time series data, a key challenge is to handle $\textit{temporal}$ dependence, in addition to cross sectional dependence amongst predictors.

Hence, techniques used in the iid case require careful extension, since in addition to dependency within the rows of the design matrix (cross sectional dependence), there is dependence withing the columns (temporal dependence). 
%Our paper seeks to bridge the gap between the classical and more general setup, where we allow dependency between the rows of the design matrix (temporal dependence) as well as that between the entries of each row (cross sectional dependence). For a comprehensive overview on recent contributions in this direction, see \citep[section 1.1]{wong2020lasso}.
%Moreover, in order to induce a general sparsity pattern, general penalty norms are extensively studied. For example, works like \citep{bach2011optimization,van2014weakly, negahban2012unified} consider norms that are $\textit{weakly decomposable}$, while work by \citep{banerjee2015estimation} considers $\textit{atomic}$ norms (these are defined in section \eqref{sec-model-and-background}). These works most commonly assume correlated designs, but they only look at correlated \textit{columns} of the design matrix, whereas we focus more on correlated $\textit{rows}$, since, in case of time series, the rows correspond to data at different time points, and are temporally dependent.
To that end,  \citet{basu2015regularized} has studies properties of LASSO regularized estimators under temporal dependence, assuming that the error term is Gaussian. \cite{wong2020lasso,wu2016performance} have studied LASSO estimates under heavier tails, and more general "mixing" conditions. 
For more general sparse structures, consistency of the $\textit{Dantzig}$ estimator is proved under Gaussianity, in \citet{banerjee2014estimation, melnyk2016estimating}, using $\textit{chaining}$ techniques \citep{talagrand2006generic}, martingales and so on (also, see \citep{nicholson2017varx} for a more applied treatment). 

In contrast to existing work, we do not use chaining arguments, but only start with a single deviation bound, then extend it using simple discretization arguments. (As a byproduct, this gives a simpler, alternate proof of consistency of the $\textit{Dantzig}$ estimator with general penalties, under Gaussianity, which may be of independent interest). Moreover, we show that the same arguments can be used to derive consistency results under heavy tails and more general notions of temporal dependence (Section \eqref{normandtempdep}). We restrict the exposition to the following two cases: $(i)$ Subweibull and $(ii)$ polynomial tails (see section \eqref{sec-Heavy-tails} for definitions). We give examples to show a strict improvement over the results in \cite{wong2016lasso, wong2020lasso, melnyk2016estimating}. In the latter case, following \cite{loh2017statistical, loh2018scale} we apply a robust regression framework with a general penalty, to derive optimal bounds. Finally, we apply these bounds to a large class of high dimensional, heavy tailed, $\textit{vector autoregressive}$ (VAR) models, under general sparsity, for which there is relatively little theoretical analysis. 

In summary, the key contributions of this work are:
\begin{itemize}
    \item We extend optimal consistency results in high-dimensional literature from the iid setup to the case of temporal dependence under a general sparsity pattern.
    \item We give novel concentration bounds \eqref{DEandRE-SW} for the $\textit{Subweibull}$ case, which can be easily compared to the iid setup.
    \item Apply our theoretical results to interesting examples of both linear high dimensional time series with SubWeibull noise.
\end{itemize}

\subsection{Organization of the paper.}

This paper is organized as follows. In section \eqref{sec-model-and-background}, we introduce the model and penalized estimation methods, as well as preliminaries related to the penalty and temporal dependence. In section \eqref{sec-Heavy-tails}, we study theoretical consistency under heavy tail--- we focus mainly on Subweibull tails, which are heavier than exponential tails but all moments still exist. The reason is, under the Subweibull assumption, the least squares loss is still fairly robust, and the temporal dependence can be made explicit in certain examples (as shown in \citep{zheng2019testing}, whose proof technique we borrow in this instance). However, if we only assume finitely many moments, then we need a robust loss function (e.g. the Huber loss) and different techniques. We also provide important applications of our results.In section \eqref{sec:SW VAR-X}, we derive optimal results for SubWeibull VAR$s$ where the transition matrix is allowed to have a general sparsity pattern. In section \eqref{sec:SW-low-rank+sparse}, we derive optimal results for SubWeibull VAR$s$ where the transition matrix is \textit{low rank}+\textit{sparse} (this requires a separate treatment due to the unique low-dimensional structure).\begin{comment}
In section \eqref{sec:SW-testing}, we give a brief overview of how our results can be used for hypothesis testing and confidence intervals.
\end{comment} 
\subsection{Notation}

Throughout the paper, we use the following notation: $||\cdot||$ denotes the $\ell_2$-norm of a vector, while $||\cdot||_p$, $||\cdot||_F$ and $||\cdot||_{nuc}$ denote the matrix norm induced by the $\ell_p$ norm $(1\le p\le\infty)$, the Frobenius  and nuclear norms of a matrix, respectively. The sparsity (or number of non-zero entries) of a matrix is denoted by $||\cdot||_0$, while$||\cdot||_{max}$ denotes the maximum of absolute values of entries of a matrix. For a $p\times d$ matrix A and $G\subseteq \{1,...,p\}\times \{1,...,d\}$ we write the submatrix $A_G = [A_{i,j} : (i,j) \in G]$. For a partition of the set $\{1,...,p\}\times \{1,...,d\}$ into disjoint groups $G_1,...,G_M$, we denote the group norm of a matrix $A$ as $||A||_{2,1} = \sum_{i=1}^M ||A_{G_i}||_F$. We also denote the mixed norm $||A||_{2,\infty}=\underset{1\le j\le d}{max} ||A_{.j}||$, where $A_{i.}$ and $A_{.j}$ denote the $i^{th}$ row and $j^{th}$ column of $A$, respectively. Let $\mbB_\mcR (0,1)=\{v:\mcR(v)\le 1\}$ and write $\mbB_2$ for $\mbB_\mcR (0,1)$ when $\mcR=||\cdot||$.  The dual of norm $\mcR()$ will be denoted by $\mcR^*()$. The cardinality of a set $J$ is denoted by $card(J)$, while its closure and convex hull is denoted by $cl(J)$ and $conv(J)$, respectively. We use $\{e_1 , \dots e_p \}$ to denote the standard canonical vectors in $\mbR^p$ with respect to the $\ell_2$ norm. For positive real numbers $a,b$, we write $b \succsim a$ if there is a positive constant $c$, independent of the model parameters, such that $b\ge ca$. Also, we say $a\asymp b$ if $a \succsim b$ and $b \succsim a$. Absolute model-free positive constants are usually denoted by $c_i$ and may change from line to line throughout the paper (these are of little concern, since they do not impact the results). The conjugate of a complex matrix $A$ is denoted as $A^*$, and if $A$ is a real matrix, its transpose is written as $A^\top$. The maximum and minimum eigenvalues of a matrix $A$ are written as $\Lambda_{\max}(A)$ and $\Lambda_{\min}(A)$. The trace and determinant of a square matrix $A$ is $\tr(A)$ and $\det(A)$, respectively.

\section{Model Formulation.}
\label{sec-model-and-background}

We start by considering the stochastic regression model  \citep{hamilton2020time}, given by  
\begin{gather}\label{stoch-reg}
y_{t}= x_t^\top \beta^* + \epsilon_{t}, \qquad t = 1,...,n,
\end{gather}
where $\{x_t,t\in\mbZ\}$ and $\{\epsilon_t,t\in\mbZ\}$ are stationary, centered linear processes possessing heavy tails, with the requirement that $x_t$ and $\epsilon_t$ are independent (expressed as $x_t\perp\epsilon_t$) for each $t$. This restriction is not very stringent and is satisfied by VAR models, which is the main focus of the paper. The covariate process $\{x_t\}$ is $p$-dimensional, while the noise process $\{\epsilon_t\}$ is univariate. The true regression parameter $\beta^*\in\mbR^p$ is unknown and needs to be estimated, based on $n$ data pairs $\{(x_t,y_t),1\le t\le n\}$. Denote by  $Y^\top=[y_1\dots y_n],X^\top=[x_1 \cdots x_n]$ and $\mcE^\top=[\epsilon_{1} \cdots \epsilon_{n}] $; then, the sample version of the model can be expressed in matrix form as 
\begin{gather} \label{stoch-reg-matrix-form}
Y = X \beta^*+ \mcE .
\end{gather}
Our focus concerns the high dimensional regime, wherein $n \ll p$. The main application of this setup will be a $\textit{vector}$ $\textit{autoregressive}$ model of lag $d$, in $\mbR^p$.  Formally, a VAR($d$) process ($d\ge 1$ represents a $\textit{fixed}$ lag), as follows: 
\begin{gather}
Z_t=B_1^\top Z_{t-1}+\dots +B_d^\top Z_{t-d} +\epsilon_t ,
\end{gather}
wherein each $B_k, k=1,\dots,d$ is a fixed transition coefficient matrix in $\mbR^{p\times p}$ and innovations $\{\epsilon_t\}$ are $p$-dimensional random vectors, whose components $\{\epsilon_{tj}\}$ are independent Subweibull($\gamma_2$) random variables with $||\epsilon_{tj}||_{\psi_{\gamma}}\le K$ for all $t$ and $j$ (see section \eqref{sec-Heavy-tails} for a definition of Subweibull norm). 

We consider the following scenarios: 
\begin{enumerate}
    \item 
    $\gamma_2=2$ (Gaussian/sub-Gaussian tails), and
    \item
    $\gamma_2\in(0,1]$ ($\textit{semi-exponential}$ tails). 
\end{enumerate}
The parameter of interest is the matrix $B^\top=[B_1^\top,\dots,B_d^\top]$ with $dp^2$ entries. Every VAR(d) process has an equivalent VAR(1) representation as $\tilde{Z}_t =\tilde{B}^\top\tilde{Z}_{t-1} +\tilde{\epsilon}_t$ where
\begin{gather}\label{eq:VAR1}
\tilde{Z_t}:=
\begin{bmatrix}
Z_t \\ 
Z_{t-1} \\
\vdots \\
Z_{t-d+1}
\end{bmatrix}_{(pd \times 1)} ,
\quad
\tilde{\epsilon}_t:=
\begin{bmatrix}
\epsilon_t \\ 
\bs{0} \\
\vdots \\
\bs{0} 
\end{bmatrix}_{(pd \times 1)}  ,
\quad
\tilde{B}^\top:=
\begin{bmatrix}
B_1^\top & B_2^\top & \cdots & B_{d-1}^\top & B_d^\top \\ 
I_p & \bs{0} & \cdots & \bs{0}  & \bs{0} \\ 
\bs{0} & I_p &        & \bs{0}  & \bs{0} \\ 
\vdots &     & \ddots &  \vdots & \vdots\\ 
\bs{0} & \bs{0} & \cdots & I_p & \bs{0} 
\end{bmatrix}_{(dp\times dp)} .
\end{gather}
Suppose one observes $T\ge d$ data points $\{Z_0,\dots,Z_T\}$; then, the original VAR(d) process can be expressed as $Z_t=B^\top\tilde{Z}_{t-1}+\epsilon_t,d\le t\le T$. Splitting this up into $p$ parallel regressions, the $j^{th}$ regression is given by
\begin{gather}
Z_{tj}= B_{j:}^\top\tilde{Z}_{t-1}+\epsilon_{tj}, \quad d\le t\le T,
\end{gather}
wherein $ B_{j:}$ represents the $j^{th}$ row of $B^\top$. Fix $1\le j\le p$ and rewrite $Z_{tj}=y_{t-d+1}$, $B_{j:}=\beta^*$, $\tilde{Z}_{t-1}=x_{t-d+1}$ and $\epsilon_{tj}=\eta_{t-d+1}$ (we suppress dependence on $j$ for now). Hence, the $j^{th}$ regression is just a stochastic regression $y_i=x_i^\top\beta^* +\eta_i,1\le i\le n$, where $n=T-d+1$, or 
\begin{gather}\label{VAR-stoch-reg}
Y=X\beta^*+\eta    
\end{gather}
in matrix form, and the parameter $\beta^*\in\mbR^{pd}$. Hence, this falls under a general stochastic regression framework (\eqref{stoch-reg}), with Subweibull tails. One of the key feature of this model is that, it is an \textit{endogenous} process---the ``input'' $\{x_t\}$ and output $\{y_t\}$ are both driven by the noise $\{\eta_t\}$. As a result the following occur: $(i)$ if the noise $\{\eta_t\}$ is not Gaussian and/or heavy tailed, this causes both $\{x_t\}$ and $\{y_t\}$ to be non-Gaussian and/or heavy tailed, respectively, and $(ii)$ the temporal dependence in the ``predictors'' $\{x_t\}$ is due to noise $\{\eta_t\}$ and further, is a function of the regression parameter $\beta^*$. Estimation of \eqref{stoch-reg} is feasible under the assumption that $\beta^*$ has a low dimensional structure; for example, it is sparse/group sparse. In that case, a regularized estimator, such as the so-called Dantzig selector \citep{candes2007dantzig}:
\begin{gather}\label{dantzig-selector-def}
\hat{\beta}= \underset{\beta \in \mbR^p }{\argmin}  \quad \mcR(\beta) \quad \text{s.t} \quad \mcR^*\left(\frac{X^\top(Y -X\beta)}{n}\right) \le \lambda_n
\end{gather}
or the LASSO \citep{tibshirani1996regression}:
\begin{gather}\label{lasso-def}
\hat{\beta} = \underset{\beta\in\mbR^p}{
\mathrm{\argmin}} \quad \frac{1}{n} ||Y - X\beta||^2 + \lambda_n \mcR(\beta) 
\end{gather}
will yield an estimate for the regression parameter, 
where $\mcR()$ is the regularizer, and $\lambda_n$ is a tuning parameter. We assume $\mcR$ is a norm, and $\mcR^*$ is its dual norm. These two estimators are essentially equivalent, if $\mcR$ corresponds to the $\ell_1$ norm \citep{bickel2009simultaneous}, or generally being $\textit{decomposable}$ \citep{wainwright2019high}. 

Hence, we can carry out $p$ regressions $\textit{simultaneously}$ with the $\textit{same}$ norm $\mcR$ and tuning parameter $\lambda_n$. The only difference will be the error probability (regret bound) for simultaneous estimation: if the regret bound for each regression is $\mfp_{err}$ (which will be the same for every regression as it only depends on $\mcR$), then the consistency result holds with probability at least $p(1-\mfp_{err})-(p-1)=1- p\mfp_{err}$ (through a simple Bonferroni bound). However, the dependence factor (which will depend on $\tilde{B}^\top$ and the Subweibull index $\gamma_2$) cannot be derived $\textit{explicitly}$ from this general framework. This is somewhat disappointing, as it is not clear $\textit{how}$ this factor is directly affected by the  target parameter $B^\top$ (due to endogeneity). Thus, in order to get an explicit expression of the dependence factor, we define a new measure of dependence in section \eqref{SubweibullVAR}, and derive consistency rates and optimal sample size, featuring this dependence factor.

\section{Background on Norms and Dependence.}\label{sec-Heavy-tails}

\subsection{Subweibull Norms.}

We start with the notion of an $\textit{Orlicz}$ norm (\citep{van2013bernstein}), that generalizes the tail decay of a random variable $Z$.

\begin{definition}

Let $\Psi:[0,\infty)\rightarrow[0,\infty)$ be an increasing and convex function with $\Psi(0)=0$. The $\Psi$-Orlicz norm of $Z$ is
\begin{gather}\label{orlicz-norm}
||Z||_\Psi:=\inf\{c>0:\mbE\parenc{\Psi(|Z|/c)\le 1}\} . 
\end{gather}
\end{definition}

Well known cases include $\textit{polynomial tails}$ that corresponds t $\Psi(z)=z^m$, $m\in\mbN$, sub-Gaussian random variables given by $\Psi(z)= \exp(z^2)-1$ and sub-exponential random variables when $\Psi(z)= \exp(z)-1$. To generalize sub-Gaussian and sub-exponential norms, let us consider $\Psi(z)= \exp(z^\gamma)-1$, for any $\gamma\in(0,2]$ (called the $\textit{tail}$ $\textit{index}$). This is called the Subweibull norm (technically a $\textit{quasi-norm}$ for $\gamma\in(0,1)$), and random variables $Z$ with finite Subweibull norm exhibit heavier tails than both sub-exponential and sub-Gaussian distributions. We shall rename the norm as $||\cdot||_{\psi_{\gamma}}$ in this case, and refer to any random variable $Z$ with finite Subweibull norm with tail index $\gamma$, as Subweibull($\gamma$). Several equivalent characterizations and properties of Subweibull random variables exist in the literature (e.g. Lemma 5, \citep{wong2020lasso}, Appendix A, \citep{gotze2019concentration}). We can then extend the notion of a Subweibull random variable to a Subweibull random $\textit{vector}$ as follows.

\begin{definition}

A random vector $X\in \mbR^p$ is said to be Subweibull($\gamma$) if $v^\top X$ is Subweibull($\gamma$) for all $v\in\mbB_2$, and its Subweibull norm is given by
\begin{gather}\label{subweibull-norm}
||X||_{\psi_\gamma}:= \underset{v\in\mbB_2}{\sup} ||v^\top X||_{\psi_\gamma} .
\end{gather}
\end{definition}

We will focus mainly on VAR models generated by a noise process with Subweibull tails, since its temporal dependence (see (\eqref{cross-corr-meas})) is easy to quantify and the obtained rates can be compared and contrasted to the Gaussian case, by simply plugging the Subweibull tail-index $\gamma=2$. For more general results on regularized stochastic regression with Subweibull noise, under $\textit{mixing}$, see Appendix \eqref{Stoch-reg-SW}.

\subsection{Stationarity}

\begin{definition}

A process $\{x_t\}$ is $\textit{strictly}$ stationary if for all $m,n,r\in\mathbb{N}$, the vector $(x_m,...,x_{m+n})$ has the same distribution as $(x_{m+r},...,x_{m+n+r})$. It is $\textit{weakly}$ (or covariance) stationary if the autocorrelations $Corr(x_t,x_{t+l})$ does not depend on $t$ for all $l\in\mbZ$. For a Gaussian process, the two notions of coincide. However, this fails to hold in general.

\end{definition}

\subsection{Role of temporal dependence: A Comparative Overview.}\label{SubweibullVAR}

From a technical standpoint, temporal dependence factors arise while using concentration inequalities for dependent data, in the course of proving consistency of our penalized estimates. Our goal is to recover optimal consistency rates that have already been derived in the independence setting, modulo a multiplicative dependence factor. Ideally, this dependence factor should be easy to interpret, and reduce to 1 for iid data (thereby being a true extension of the latter). Unfortunately, there is no single, unified framework for temporal dependence that gives us tight concentration bounds for every family of time series. Hence, we require different, but related, notions of dependence that is suited to specific examples. These different notions of dependence do not imply one another, and may hence be seen as complementary to each other. 

From the standpoint of consistency rates, the temporal dependence factor ``inflates" the consistency bounds, and the finite, minimum sample size required to achieve the bound that holds for the independent case. Hence, the effect of dependence on the penalized estimates is clear.

\subsubsection{Temporal Dependence Measure for Sub-Weibull Linear Processes.}\label{temp-dep}

We first define a measure of temporal dependence for a Subweibull linear process, that will be used when deriving subsequent consistency results. Consider a linear filter of an innovation process $\{\eta_t\}\sim\text{IID}(0,\Sigma_\eta)$, given by $x_t=\sum_{i\ge 0}A_i\eta_{t-i} =\mcA(\mathfrak{B})\eta_t$ where $\mathfrak{B}$ denotes the backshift operator, the matrix power series $\mcA(z)=\sum_{i\ge 0}A_i z^i$ defined on the complex plane $z\in\mbC$, satisfy $\det(\mcA(z))\neq 0$ on the unit disk $\{z\in\mbC:|z|\le 1\}$  (stationarity). We then define the following and assume it is finite:
\begin{gather}\label{cross-corr-meas}
\msfC(\mcA):=\sum_{i\ge 0}\sum_{j\ge 0}||A_{i+j}||_2 ||A_j||_2 . 
\end{gather}
Note that $\msfC(\mcA)<\infty$ implies $\sum_{i\ge 0}||A_i||<\infty$, i.e. the process is stable. This definition is motivated by the fact that the dependence in the linear process $\{x_t\}$ is due to the auto-correlation terms, i.e. $\cov(x_0,x_l) =\Sigma_X(l)=\sum_{i=0}^\infty A_i\Sigma_\eta A_{i+l}^\top$ for all non-negative lags $l$. For a connection between this measure of dependence, and the one introduced in \citep{basu2015regularized}, see Section \eqref{normandtempdep}. In particular, for a stationary VAR(1) process $x_t=Ax_{t-1}+\eta_t$, we have the causal representation $x_t=\sum_{i\ge 0}A^i\eta_{t-i}$. Assuming the spectral radius of the transition matrix $\rho(A)<1$ (stability, see \citep{lutkepohl2005new}), it is easy to verify that  the series $\sum_{i\ge 0}||A^i||_2<\infty$ (using Gelfand's formula: $\rho(A)= \underset{k\rightarrow\infty}{\lim} ||A^k||_2^{1/k}$). Then,
\begin{align}
\msfC(\mcA)=\sum_{i\ge 0}\sum_{j\ge 0} ||A^{i+j}||_2 ||A^j||_2 
& \le\sum_{i\ge 0}\sum_{j\ge 0} ||A^i||_2 ||A^j||_2^2 \\
& =\sum_{i\ge 0} ||A^i||_2 \sum_{j\ge 0}||A^j||_2^2 
\le\left[\sum_{i\ge 0} ||A^i||_2 \right]^3<\infty .
\end{align}
We stress that the assumption: $\rho(A)<1$ is mild and standard. Many authors (e.g. \citep{loh2012high,bickel2008covariance}) work with the much stronger assumption that the spectral norm $||A||_2<1$. However, this fails for VAR(d) models of lag $d>1$ (\citep*[lemma E.1]{basu2015regularized}), and therefore can not be applied in general. Even with a VAR(1) process, we can take the transition matrix $A=
\begin{bmatrix}
a,b\\
0,a
\end{bmatrix}$,
where $|a|<1,b\in\mbR$. The spectral radius is $\rho(A)=|a|<1$, while the spectral norm is $||A||_2^2=a^2+\frac{b^2}{2}+\frac{|b|}{2}\sqrt{4a^2+b^2}\rightarrow\infty$ as $b\rightarrow\infty$, with $a$ fixed. However, our dependence measure does not have this drawback. If we do have $||A||_2<1$, then $\msfC(\mcA)\le(1-||A||_2)^{-2}$. Also, as the spectral radius $\rho(A)\rightarrow 1$, the process becomes highly unstable and the dependence factor $\msfC(A)$ blows up. For a simple VAR(1) example, if  we take $A= \begin{bmatrix}
\rho, 0\\
0, \rho
\end{bmatrix}$
where $|\rho|<1$, then, the spectral radius of A is $\rho(A)=\rho$ while the dependence factor $\msfC(A)=(1-\rho)^{-1}(1-\rho^2)^{-1} \rightarrow\infty$, as the spectral radius $\rho\rightarrow 1$.
Finally, we note that, in case of a VAR(d) process $Z_t=B_1^\top Z_{t-1}+\dots +B_d^\top Z_{t-d} +\epsilon_t$, we can write it as a VAR(1) process $\tilde{Z}_t= \tilde{B}^\top\tilde{Z}_{t-1} +\tilde{\epsilon}_t$, as in section \eqref{sec-model-and-background}. Stability of $\{Z_t\}$ implies stability of $\{\tilde{Z}_t\}$ (see e.g. \citep*[Ch 2.1]{lutkepohl2005new}), which means $\msfC(\tilde{\mcB})<\infty$ (note that $\msfC(\mcB)$ and $\msfC(\tilde{\mcB})$ are not the same). From section \eqref{sec-model-and-background}, breaking up this model into component regressions and considering a single component regression $y=X\beta^*+\eta$, we give novel concentration bounds for the deviation term $X^\top\eta/n$ and the sample Gram matrix $X^\top X/n$, which serve as starting points for the general deviation and Restricted Eigenvalue conditions (similar e.g., to Proposition \eqref{ConcBasu} in \cite{basu2015regularized} for the strictly Gaussian case).

\begin{proposition}\label{DEandRE-SW}

%\textcolor{red}{ADD A FULL STATEMENT. E.g. SPECIFY WHERE $X$ COMES FROM AND WHAT ITS PROBABILISTC FRAMEWORK IS, e.g., SUBWEIBULL}

\noindent
Consider the stochastic regression \eqref{VAR-stoch-reg} with Subweibull tails. The Gram matrix $X^\top X/n$ and the deviation term $X^\top\eta/n$ obtained from the posited model satisfy the following, respectively:
Fix $u\in\mbR^{dp}$ with $||u||\le 1$. Then, for $t>0$,
\begin{gather}
\mbP\parenf{|u^\top(X^\top X/n-\Sigma_X)u|> K^2\msfC(\tilde{\mcB})t} 
\le 6\exp\parens{-c\min\parenc{(nt)^{\frac{\gamma_2}{2}}, nt^2)}} \\
\mbP\parenf{|u^\top X^\top\eta/n|>K^2\msfC(\tilde{\mcB})t} 
\le 6dp\exp\parens{-c\min\parenc{(nt)^{ \frac{\gamma_2}{2}}, nt^2)}}.
\end{gather}
\end{proposition}

\begin{remark}
When $\gamma_2=2$, Proposition \eqref{DEandRE-SW} matches the deviation bounds \eqref{ConcBasu} for the strictly Gaussian case in \citep{basu2015regularized}, and moreover extends it to the Sub-gaussian case also (note that the simple trick of rotating a Gaussian vector to make its components independent, will not work outside Gaussianity). This closely mimics the bounds derived for the Sub-Gaussian case in \citep{zheng2019testing} and extends those bounds to heavy tails.
\end{remark}

\begin{remark}
Unlike Proposition \eqref{ConcWong} in \cite{wong2020lasso} for the more general case that requires $\textit{mixing}$ conditions, these concentration bounds hold for all $t>0$, $n\ge 1$, so no preconditions are required. Further, they are used in the sequel to obtain clean expressions for the sample size $n$ requirement, tuning parameter $\lambda_n$, as well as temporal dependence, that directly match those derived in \citep{basu2015regularized} for the strictly Gaussian case.

\end{remark}

\begin{remark}
These results are established by essentially a "truncation" argument analogous to that used for proving \citep[Lemma 5.2]{zheng2019testing}. This is possible because of the highly specific dependence structure of VAR --- it is a causal (infinite) linear combination of \textit{independent} noise/shocks. Hence, the trick is to truncate the linear series at some \textit{finite}  time point in the past and use a Hanson-Wright type inequality \citep{vershynin2018high} for this \textit{finite} linear combination of \textit{independent} noises. Finally, what's "left-over" is controlled suitably by choosing this finite time-point carefully.
\end{remark}

%\textcolor{red}{DO WE WANT TO EMPHASIZE ANYTHING ON HOW THE RESULTS ARE ESTABLISHED? TECHNICAL INNOVATION?}

\begin{remark}
We have also compared this concentration bound to those found in the most relevant literature, namely \citep{basu2015regularized} and \citep{wong2020lasso}, as well as the different quantifying dependence factors, in Appendix \eqref{sec-Compare-Dependence}.
\end{remark}

\subsection{Concepts Related to the Penalty Term.}\label{normandtempdep}

As noted earlier, the focus in the literature has primarily been on sparse/group sparse penalties for the regression coefficient $\beta^*$ in \eqref{VAR-stoch-reg} and the transition matrix $\tilde{B}$ in \eqref{eq:VAR1}. Here, we review some concepts that would be used in the sequel to establish results for much more general penalty terms that are useful in practical settings.

\begin{definition}\label{GaussianWidth}
Let G be a $p \times d$ random matrix with iid $N(0,1)$ entries. For a a set $\msT\subseteq\mbR^{p\times d}$, the Gaussian width of $\msT$ is defined as
\begin{gather}\label{gaussian-width}
w(\msT)= \underset{W\in\msT}{\sup}\tr(W^\top G) .   
\end{gather}
It measures the size of a (usually convex) set in the Euclidean space. A key challenge is to evaluate Gaussian widths of sets related to the regularizer $\mcR$.
\end{definition}

\begin{definition}\label{CompatibilityConst}
Given a set $\msC\subseteq\mbR^{p\times d}$, and a generic norm $\mcR()$ on $\mbR^{p\times d}$, the subspace norm compatibility constant is given by 
\begin{gather}\label{compatibility-constant}
    \Phi_\mcR(\msC)= \underset{W \in\msC-\{0\}}{\sup} \frac{\mcR(W)}{||W||_F}.
\end{gather}
When $d=1$, the definition reduces to the one given in \citep{negahban2012unified}. It measures the  relative price paid for switching between a generic norm $\mcR()$ and the usual Euclidean norm $||\cdot||_F$. Also, the reverse norm compatibility is given by
\begin{gather}\label{reverse-comp}
    \bar{\Phi}_\mcR (\msC)= \underset{W \in\msC-\{0\}}{\sup} \frac{||W||_F}{\mcR(W)} .
\end{gather}
In particular when $\msC=\mbB_\mcR(0,1)$, we write $\bar{\Phi}_\mcR(\msC)$ as simply $\bar{\Phi}_\mcR$. For example, we note that $||v||\le\bar{\Phi}_\mcR \mcR(v)\le \bar{\Phi}_\mcR$ for all $v\in \mbB_\mcR(0,1)$. We further assume that $\bar{\Phi}_\mcR$ is bounded above by an absolute constant (this holds in most cases under consideration).
\end{definition}

%\textcolor{red}{THE MATERIAL IN SECTIONS 3.4 and 3.6 IS INTERESTING, BUT NOTE CORE. ADD HERE A CLARIFYING REMARK AND MOVE THESE SECTIONS TO AN APPENDIX.}

\section{Subweibull VAR under general sparsity.}
\label{SubweibullVARex}

In this section, we leverage Proposition \eqref{DEandRE-SW} to obtain optimal consistency rates for different examples of sparse, SubWeibull VAR$s$. For model \eqref{VAR-stoch-reg}, we consider a general sparse structure on $\beta^*$, induced through the Dantzig selector \eqref{dantzig-selector-def} with an appropriate penalty $\mcR()$. To consistently estimate $\beta^*$, it is standard practice in high-dimensional literature (e.g.
\citep{van2011adaptive,bickel2009simultaneous}) to verify a first order "deviation" condition, and a second order "Restricted Eigenvalue (RE)" condition. The deviation condition essentially restricts the cross product term $\frac{X^\top\eta}{n}$ around zero with high probability, while the RE condition ensures the sample gram matrix $\frac{X^\top X}{n}$ is uniformly bounded away from zero over a small cone, with high probability. We prove both under the general penalty $\mcR()$.

\begin{proposition}\label{DeviationSW-VAR} {Deviation Condition for Subweibull VAR.}
Consider model \eqref{VAR-stoch-reg} posited in section \eqref{sec-model-and-background}. Then, there exists an absolute constant $c>0$ such that, for $n\ge n_{dev}:= \parens{c w^2(\mbB_\mcR(0,1))}^{4/\gamma_2-1}$, we get
\begin{gather}
\mbP\left[ \mcR^*\parenf{\frac{X^\top\eta}{n}}\ge 2\bar{\Phi}_\mcR K^2\msfC(\tilde{\mcB})\sqrt{\frac{cw^2(\mbB_\mcR (0,1))}{n}} \right]
\le 6\exp[-w^2(\mbB_\mcR(0,1))+\log dp] .
\end{gather} 
\end{proposition}

\begin{proposition}\label{RE-SW-VAR}{RE condition for Subweibull VAR.}
Consider model \eqref{VAR-stoch-reg} posited in section \eqref{sec-model-and-background}. Assume $\Lambda_{\min}(\Sigma_X)>0$.  Then, there is an absolute constant $c'>0$ such that, for a minimum sample size of 
\begin{gather}
n\ge n_{RE}:=\parenf{c'\max\parenc{1,\frac{16\bar{\Phi}_{\mcR}^2 K^4\msfC^2(\tilde{\mcB})}{\Lambda_{\min}^2(\Sigma_X)}} w^2[\mcT\cap\mbB_2)}^{2/\gamma_2}    
\end{gather}
we obtain
\begin{gather}
\mbP\parens{\underset{v\in\mcT\cap\mbB_2}{\inf}\frac{v^\top X^\top X v}{n}\ge\alpha_{RE}}\ge 1-6\exp[-w^2(\mcT\cap\mbB_2)],
\end{gather}
where the restricted eigenvalue is $\alpha_{RE} =\Lambda_{\min}(\Sigma_x)/2$. 

\end{proposition}

Using these deviation and RE conditions, we can consistently estimate each of the $p$ regressions that the Subweibull VAR model can be decomposed into. Subsequently, we can combine these estimates using a Bonferroni bound that leads to the following optimal consistency result.

\begin{proposition}\label{VAR-Subweibull-consistency}

Consider the VAR(d) model posited in section \eqref{sec-model-and-background}. Further, suppose the minimum sample size and tuning parameter satisfies
\begin{gather}
n\ge\max\parenc{n_{dev},n_{RE}}, \quad
\lambda_n =2\bar{\Phi}_\mcR K^2\msfC(\tilde{\mcB})\sqrt{\frac{cw^2(\mbB_\mcR (0,1))}{n}} .
\end{gather}
Then, denoting the columns of $B$ (equivalently rows of $B^\top$) as $B_1,\dots,B_p$, and assuming they share a common sparsity pattern, the augmented penalized estimate $\hat{B}$ satisfies
\begin{gather}
\underset{1\le j\le p}{\max}||\hat{B}_j-B_j||\leq \frac{2\lambda_n\Phi_\mcR(\mcT)}{\alpha_{RE}}, \\
\underset{1\le j\le p}{\max} \mcR(\hat{B}_j-B_j) \le\frac{2\lambda_n\Phi_\mcR^2(\mcT)}{\alpha_{RE}},
\quad\text{(Estimation error)}, \\
\underset{1\le j\le p}{\max} (\hat{B}_j-B_j)^\top \frac{X^\top X}{n}(\hat{B}_j-B_j)\le \frac{4\lambda_n^2\Phi_\mcR(\mcT)}{\alpha_{RE}} 
\quad\text{(Prediction error)}.
\end{gather} 
with probability at least $1-6dp^2\exp[-w^2(\mbB_\mcR(0,1))]-6p\exp[-w^2(\mcT\cap\mbB_2)]$. 

Noting that $x_i=\tilde{Z}_{i+d-2}$, for $1\le i\le n$, where $n=T-d+1$ is the sample size, the design matrix $X$ is given by
\begin{gather}
X =
\begin{bmatrix}
\tilde{Z}_{d-1}^\top \\ 
\tilde{Z}_{d-2}^\top \\ 
\vdots\\ 
\tilde{Z}_{n+d-2}^\top 
\end{bmatrix}_{(dp\times dp)} = \quad
\begin{bmatrix}
Z_{d-1}^\top & Z_{d-2}^\top & \cdots & Z_0^\top \\ 
Z_d^\top    & Z_{d-1}^\top & \cdots & Z_1^\top \\ 
\vdots &     \vdots         & \ddots &   \vdots \\ 
Z_{T-1}^\top & Z_{T-2}^\top & \cdots & Z_{T-d}^\top \\  
\end{bmatrix}_{(dp\times dp)} .
\end{gather}
The restricted eigenvalue is given by $\alpha_{RE}=\Lambda_{min}(\Sigma_X(0))/2=\Lambda_{min} (\Sigma_{\tilde{Z}})/2$.

\end{proposition}

The key algebraic quantities of interest are the following:

\begin{itemize}
    \item The Gaussian width of the unit norm ball : $w(\mbB_\mcR(0,1))$.
    \item The Gaussian width of the spherical cap of the tangent cone $\mcT$ : $w(\mcT\cap\mbB_2)$.
    \item The subspace compatibility constant $\Phi_\mcR(\mcT)$. 
    \item The reverse compatibility constant : $\bar{\Phi}$. 
\end{itemize}

Next, we provide estimates for these quantities for different examples of the regularizer $\mcR$.

\begin{example}\label{Dantzig-l1norm}{that  considers the $\ell_1$ norm and serves for illustration purposes.}

Suppose the parameter $\beta^*$ is $s$-sparse, i.e. $||\beta^*||_0=s$. Then, the $\ell_1$ norm is a convex relaxation of $||\cdot||_0$. Hence, for $\mcR()=||\cdot||_1$, (see Examples 1.1, 2.1, 3.1, \citep*{banerjee2015estimation} and Proposition 3.10, \citep*{chandrasekaran2012convex}), we have 
\begin{gather}
w(\mbB_\mcR(0,1))\le 2\sqrt{\log 2p},\quad
w^2(\mcT\cap\mbB_2)\le 2s\log(p/s) + \frac{5}{4}s,\quad\\
\Phi_\mcR(\mcT)\le 2\sqrt{s},\quad
\bar{\Phi}\le 1.
\end{gather}
Hence, ignoring dependence factors, the $\ell_2$ consistency rate is $\mcO(\sqrt{s\log p/n})$, and the minimum sample size requirement is $n\succsim\max \parenc{s(\log(p/s)+1),\log 2p}$ (which matches Proposition 3.3, \citep*{basu2015regularized}, which was derived for the $\ell_1$-LASSO).

\end{example}

\begin{example}\label{Dantzig-OWLnorm}{Sorted $\ell_1$ norm.}

The $\ell_1$ norm has the drawback that it treats all  coefficients in the $\beta$ vector equally; this is a problem if there is there is cross-sectional dependence among the significant components of $\beta^*$. In that case, the $\ell_1$ penalty basically selects an arbitrary subset of the significant components, whereas it is desirable to include \textit{all}  relevant variables in the analysis. Several authors have addressed this problem; the most popular solution being the elastic net. Another solution was proposed by \citep*{bondell2008simultaneous} to deal with sparse regression with correlated variables, and generalized by \citep*{bogdan2013statistical} to a more general class of penalties characterized by the Ordered Weighted $\ell_1$ or OWL norm (also called the SLOPE, e.g. see \citep*{bogdan2015slope,stucky2018asymptotic}). It is defined as
\begin{gather}
\mcR(\beta)=\sum_{i=1}^p \mfw_i\beta_i^{\downarrow}  
\end{gather}
for some weights $\mfw_1\ge...\ge\mfw_p\ge 0$, wherein  $\{\beta_i^{\downarrow},1\le i\le p\}$ is a decreasing arrangement of $\{|\beta_i|,1\le i\le p\}$. Even though this norm is \textit{not decomposable}, it is $\textit{atomic}$ \citep*{zeng2014ordered}, and helps cluster significant variables when they exhibit strong cross-sectional dependence \citep*{figueiredo2014sparse, figueiredo2016ordered}. Let $\bar{\mfw}_i=(\sum_{t=1}^i \mfw_t)/i$. Then $\mfw_1\ge...\ge\mfw_p\ge 0$ implies $\bar{\mfw}_1\ge...\ge\bar{\mfw}_p\ge 0$. Denote $\tilde{\mfw}_s$ as the average of $\mfw_{s+1},...,\mfw_p$. Also, $\mcR(v)\ge \mfw_1||v||_1\ge\mfw_1 ||v||$.  Then, for $s$-sparse  $\beta^*$ we get (see \citep*{banerjee2015estimation})
\begin{gather}
w(\mbB_\mcR(0,1))\le 2\underset{1\le i\le p}{\max} \parenc{\frac{\sqrt{2+\log(2p/i)}}{\bar{\mfw}_i}} + 2\sqrt{\log 2p}/\bar{\mfw}_p 
\le 2\sqrt{2+\log 2p}/\bar{\mfw}_p, \\
w^2(\mcT\cap\mbB_2)\le\frac{2\mfw_1^2}{\tilde{\mfw}_s} s\log(p/s) + \frac{3}{2}s,\quad\\
\Phi_\mcR(\mcT)\le \frac{2\mfw_1^2}{\tilde{\mfw}_s} \sqrt{s},\quad \bar{\Phi}\le\mfw_1^{-1}.
\end{gather}
Thus, the order of consistency is $\mcO\parenf{2\mfw_1^2/ (\bar{\mfw}_p\tilde{\mfw}_s)\sqrt{(s\log p)/n}}$. The minimum sample size required is given by $n\succsim\max \parenc{(\mfw_1^2/\tilde{\mfw}_s) s(\log(p/s)+1), \log(2p)/\bar{\mfw}_p^2}$. When $\mfw_1=...=\mfw_p=1$, we get $\bar{\mfw}_p =\tilde{\mfw}_s =1$ and $\mcR$ reduces to the $\ell_1$ norm, with the usual rate of consistency.
\end{example}

\begin{example}\label{Dantzig-group-norm}{Group sparsity.}

Suppose, instead of element-wise sparsity, we assume a group structure on the parameter $\beta^*$. If the parameter space comprises of $M$ (possibly overlapping) groups $\mfG=\{G_1,...,G_M\}$, with $G_i$ being the set of parameter indices in the $i^{th}$ group and $\cup_{i=1}^M G_i=\{1,..,p\}$, then it is appropriate to consider the $\textit{atomic}$ norm $\mcR$ induced by this grouping $\mfG$. When the groups $G_i$ are disjoint, then $\mcR$ is just the $\ell_{2,1}$ norm given by $\mcR(\beta)= ||\beta||_{2,1} =\sum_{i=1}^M||\beta_{G_i}||$. Let the maximum group size be denoted by $m$ and suppose only $s$ many groups are active. Then, by Lemma 2 in \citep*{banerjee2015estimation}, and \citep*{rao2012universal}, we obtain
\begin{gather}
w(\mbB_\mcR(0,1))\le\sqrt{m}+2\sqrt{\log M},\quad
w^2(\mcT\cap\mbB_2)\le\parenc{(\sqrt{2\log(M-s)}+ \sqrt{m})^2+m}s,\quad\\
\Phi_\mcR(\mcT)\le s,\quad\text{(overlapping groups)}\quad
\text{or}\quad \sqrt{s},\quad\text{(non-overlapping groups)} \\
\bar{\Phi}\le 1.
\end{gather}
For non-overlapping groups, we then have the usual rate $\mcO\parenf{\sqrt{s}(\sqrt{m} +\sqrt{\log M})/\sqrt{n}}$ and the sample size scales as $n\succsim s(m+\log M)$.

\end{example}

\begin{example}\label{Dantzig-ksupport-norm}{The $k$-support norm. }

The ``elastic net" regularizer comprising of a combination of the $\ell_1$ and $\ell_2$ norms \citep*{zou2005regularization} is often advocated as a better alternative to the $\ell_1$ LASSO, since it may be considered as a convex relaxation of $||\cdot||_0$ with the scale set by the $\ell_2$ norm. However, a tighter convex relaxation is possible, and captured by the so-called ``k-support'' norm introduced in \citep*{argyriou2012sparse}. It is an \textit{atomic} norm, shown to recover sparse parameters $\beta^*$ better than the elastic net empirically. We have, from \citep*{banerjee2015estimation}, Proposition 3.1, \citep*{argyriou2012sparse},
\begin{gather}
w(\mbB_\mcR(0,1))\le\sqrt{k}+2\sqrt{k\log(p/k)+k},\quad
w^2(\mcT\cap\mbB_2)\le
\sqrt{\frac{2\beta_{\max}^*}{\beta_{\min}^*}s\log(p/s) +\frac{3}{2}s} \\
\Phi_\mcR(\mcT)\le \sqrt{2}\parenf{1+ \frac{2\beta_{\max}^*}{\beta_{\min}^*}} ,\quad 
\bar{\Phi}\le 1.
\end{gather}
Here, $\beta_{\max}^*=\underset{1\le j\le p}{\max}|\beta_j^*|$ and $\beta_{\min}^*=\underset{1\le j\le p}{\min}|\beta_j^*|$.

\end{example}

\subsection{Comparisons to Previous Work.}

The result in Proposition \eqref{VAR-Subweibull-consistency} indicates how the sample size $n$ and the tuning parameter $\lambda_n$ scale by a factor $\msfC^2(\tilde{\mcB})$, the latter being \textit{computable}. For ease of presentation, suppose the lag $d=1$, so that $Z_t=B^\top Z_{t-1}+\epsilon_t$. Assuming $||B||_2<1$ (e.g. symmetric $B$), we get $\msfC(\mcA) \le(1-||B||_2)^{-2}$. Thus, as $||B||_2\rightarrow 1$, the process becomes highly unstable which causes the sample size $n$ to blow up. On the other hand, this bound is tight in the sense that, in case of independence (take $B=0$ for a white noise process), $\msfC=1$. 

\medskip
\noindent
\textbf{Comparison with \citep{basu2015regularized}.}
The special case of $\gamma_2=2$ with the penalty $\mcR$ taken to be the $\ell_1$ norm, leads to a direct comparison with that work that assumes $\textit{Gaussianity}$, which in turn simplifies calculations. Specifically, the switch from dependence $\rightarrow$ independence can be accomplished though the following: if a vector $X\sim\mcN(0,\Sigma)$, then $\Sigma^{-\frac{1}{2}}X \sim\mcN(0,\mcI)$. The price of dependence is simply $\Sigma^{-\frac{1}{2}}$ or, in terms of a number -- $\Lambda_{min}(\Sigma^{\frac{1}{2}}$ or $\Lambda_{max}{\Sigma^\frac{1}{2}}$). However, this is specific to Gaussian processes and fails even for a dependent sub-Gaussian process (e.g. when the data is discrete or have finite range). Nonetheless, we expect the optimal $\ell_1$-LASSO rate $\sqrt{\frac{s\log p}{n}}$ even in this case. We recover this bound for the sub-Gaussian case also, simply by putting $\gamma_2$=2 in (\eqref{VAR-Subweibull-consistency}). Moreover, \citep{basu2015network} considers the popular $\ell_1$ penalty which has the advantage of being \textit{decomposable}, whereas this work covers more general classes of penalties and hence the obtained result is more widely applicable.

\medskip\noindent
\textbf{Comparison with  \citep{wong2020lasso}.}
Note that the obtained result is \textit{strictly tighter} than its counterpart presented in Appendix E in the aforementioned paper for the Subweibull regression framework under $\textit{mixing}$ conditions. For direct comparison, consider an $\ell_1$ penalty. The process $\{Z_t\}$ can be shown to be geometrically $\beta$-mixing, with geometric index $\gamma_1=1$. Leveraging the mixing framework in \citep{wong2020lasso} to this example, will lead to a sub-optimal sample size requirement \citep*[Corollary 9]{wong2020lasso}. Specifically, in the sub-Gaussian case, plugging $\gamma_2=2$, $\gamma_1=1$ (and taking the $\ell_1$ penalty) in Propositions \eqref{dev-Subweibull}, \eqref{RE_condition_Subweibull} and \eqref{det-errbd-SW-stochreg}, will imply $\gamma=1/2$, leading to the inferior bound $n\succsim \max\parenc{(\log p)^3,(s\log p)^2}$ (see \citep*[ Corollary 4]{wong2016lasso}), instead of the optimal bound $n\succsim \max\parenc{\log(ep/s), s\log p}$ (ignoring dependence factors). Moreover, in the latter work, the temporal dependence factors are not explicitly computable or easily interpreted. This is due to the generality of the mixing framework, which is simultaneously its strength and weakness. Its strength lies in the fact that it is applicable to a very large family of time series; specifically, those admitting a "Wald representation". On the other hand, its generality precludes obtaining tight concentration inequalities without further assumptions. This issue is avoided in the case of a sub-Weibull linear process that enables us to obtain simpler and tighter concentration bounds with easy to interpret dependence factors.
%This problem may be unavoidable based on the cleanest results (as far as we know) on concentration bounds involving heavy tailed distributions under temporal dependence. It would be interesting to derive tight Bernstein or Hoeffding inequalities for mixing processes. However, we can circumvent the problem when we are dealing with a specific family e.g. a SubWeibull linear process--- we do not require that level of generality. Instead, we can exploit the specific properties of, e.g SubWeibull VAR$s$, to get simpler, tighter concentration bounds, with easy-to-interpret dependence factors.

\medskip\noindent
\textbf{Comparison with  \citep{melnyk2016estimating}.} 
Proposition \eqref{VAR-Subweibull-consistency} extends results on regularized VAR models assuming sub-exponential tails for the error distribution, obtained in that paper, to the sub-Weibull case. Note that \citep{melnyk2016estimating} employ $\textit{generic chaining arguments}$ \citep[Theorem 1.2.7, 1.2.9]{talagrand2006generic} in the derivations. It 
is doubtful whether the arguments can be repeated when considering \textit{beyond} sub-exponential tails - e.g $\textit{semi}$-exponential tails, i.e., $\gamma_2<1$. In contrast, our proof techniques are simpler and directly lead to tractable dependence factors.

\section{Subweibull VAR with a Low Rank plus Sparse Transition Matrix.}\label{sec:SW-low-rank+sparse}

Next, we consider, for simplicity, a VAR(1) model where the transition matrix $B$ is low-rank plus group sparse given by
\begin{gather}\label{eqn:model-sparse-lowrank}
Z_t= B^\top Z_{t-1}+\epsilon_t, \quad 
B=L^*+R^*,\quad\text{rank($L^*$)=r},
\end{gather}
where $L^*$ represents the low rank component and $R^*$ represents either a sparse $S^*$, or group-sparse component
$G^*$. We assume the number of non-zero elements in the sparse case is $\|S^*\|_0=s$, while in the group sparse case there are $g$ non-zero groups out of $M$ groups, with $r\ll p, s\ll p^2$ and $g\ll p^2$. The matrix $L^*$ captures a common dense persistence structure across \textit{all} $p$ component series, specifically, strong cross-correlations between the component series, which a simple sparse VAR model cannot handle. Also, we assume the noise components $\{\epsilon_{tj}\}$ to be independent SW($\gamma_2$) with $||\epsilon_{tj}||_{\psi_{\gamma_2}} \le K,\forall i,j$. We want to estimate $L^{*}$ and $R^{*}$ accurately based on a sample $\{Z_0,\dots,Z_T\}$, when $N \ll p^2$.
\begin{gather}\label{eqn:data-sparse-lowrank}
\underbrace{\left[\begin{array}{c}
Z_T^\top \\
\vdots \\
Z_1^\top \end{array} \right]}_{\mcY} =  
\underbrace{\left[ \begin{array}{c}
Z_{T-1}^\top\\
\vdots \\
Z_0^\top\end{array} \right]}_{\mathcal{X}} B
+ \underbrace{\left[ \begin{array}{c} \epsilon_T^\top \\
\vdots \\
\epsilon_1^\top \end{array}\right]}_{E}.
\end{gather}
This is a standard matrix regression problem with $n=T$ samples and $p^2$ parameters. However, we face an \textit{identifiability} issue in the estimation of the low rank and sparse components $L^*$ and $R^*$. For example, if the low-rank component $L^*$ itself is $s$-sparse and the sparse component $R^{*}$ is of rank $r$, then we cannot hope for a way to estimate $L^{*}$ and $R^{*}$ separately  without further restrictions. Thus, we impose the restriction that condition that the low rank part should not be too sparse and the sparse or group-sparse part should not be low-rank (e.g. \citep{chandrasekaran2011rank, agarwal2012noisy}). We thus have the following estimation procedure:
\begin{gather}
(\hat{L},\hat{R})= \argmin_{\substack{L,R\in\mbR^{p\times p}  
\\ L\in\Omega}} \frac{1}{2}\left\|\mcY-\mcX (L+R)\right\|^2_F + \lambda_n\|L\|_*+\mu_n\|R\|_{\diamond} ,
\end{gather}
where $\Omega=\{L\in\mbR^{p\times p}: \|L\|_{\max}\le \alpha/p\}$ (for sparse) or $\{L\in\mbR^{p\times p}:\|L\|_{2,\max}\le\beta/\sqrt{M}\}$ (for group sparse), $\|\cdot\|_{\diamond}$ represents $\|\cdot\|_1$ or $\|\cdot\|_{2,1}$ depending on sparsity or group sparsity of $R$, and $\lambda_n$ and $\mu_n$ are tuning parameters. The parameters $\alpha$ and $\beta$ control for the degree of ``non-identifiability'' of the matrices allowed in the model. Specifically, large $\alpha$ provide sparser estimates of $S$, while allowing simultaneous sparse and low-rank components to be absorbed in $\hat{L}$. On the other hand, smaller $\alpha$ pushes the simultaneous low-rank and sparse components to be absorbed in $\hat{S}$. The problem under independence was studied by \citep{agarwal2012noisy}, and extended to the dependent, Gaussian case by \citep{basu2019low}. We extend these results to the case of dependent data--- sub-Gaussian \textit{and} heavy tailed (tricks for Gaussianity can not be exploited). The crucial change would be an application of Proposition \eqref{DeviationSW-VAR}, to get probability bounds matching \citep*[Proposition 3]{basu2019low}. Also, following \citep{basu2019low}, we choose $\alpha$ and $\beta$ in the range $[1,p]$ and $[1,K]$, respectively. To fix ideas, let us consider low rank+sparse set up so that $R=S$.

\begin{proposition} \label{VAR-L-S}
Consider the low rank+ sparse VAR with Subweibull noise posited by \eqref{eqn:model-sparse-lowrank}.
%\textcolor{red}{SAY FOR THE MODEL POSISTED IN EQN XXX AND ASSUMING (LIKE IN THE PREVIOUS PROPOSITION)}
There are absolute constants $c_i>0$ such that for
\begin{gather}
n\ge\parens{\frac{c_0 K^4\msfC^2(\mcB)}{\Lambda_{\min}^2(\Sigma_X)} \cdot p}^{4/\gamma_2-1}, \quad \|L^*\|_{\max}\le\alpha/p,     
\end{gather}
we have
\begin{gather}
\|\hat{S}-S^*\|^2_F +\|\hat{L}-L^*\|^2_F \le 
c_1\parens{\frac{K^4\msfC^2(\mcB)}{\Lambda_{\min}^2(\Sigma_X)}
\frac{(rp +s\log p)}{n} +\frac{s\alpha^2}{p^2} }.
\end{gather}
with probability at least $1- c_2\exp[-c_3\log p]$.

\end{proposition}

The result in Proposition \eqref{VAR-L-S} is new the literature. As noted in \citep{basu2019low}, the first term in the estimation error is due to the randomness in the data and limited sample size, and becomes small as the sample size increases. The second term is due to the unidentifiability of the problem, and does not vanish, even as the sample size grows large. Of course, plugging $\gamma_2=2$, we get back the result given in \citep*[Proposition 4]{basu2019low} (note that this covers the $\textit{sub}$-Gaussian case also, which can not be handled with the same tools as the Gaussian case). The results for group-sparse component and extension to VAR models with lags bigger than 1, is also straightforward following along the lines in \citep{basu2019low}.

\section{Subweibull VAR with Exogenous Predictors (VAR-X).}\label{sec:SW VAR-X}

In many applications, a VAR’s forecasts can be improved by incorporating variables which are determined outside of the VAR. Examples of exogenous variables include leading indicators, weather-related measurements, global macroeconomic variables such as world oil prices, etc. Econometricians call these models “VAR-X”, or “transfer function” or “distributed lag” models. VAR-X finds popularity in the modeling of small open economies, as they are generally sensitive to a wide variety of global macroeconomic variables which evolve \textit{independently} of their internal indicators (e.g. \citep{cushman1997identifying}). Furthermore, VAR-X models are applied just in marketing \citep{nijs2007retail}, political science \citep{wood2009presidential}, and real estate \citep{brooks2000forecasting}. Obviously the curse of dimensionality that exists in VAR is compounded in VAR-X models due to more (this time, exogenous) variables. To this end, Nicholson et al \citep{nicholson2017varx,nicholson2020high} have given examples of various penalties $\mcR()$ under the Gaussian setup, that sparsify the problem and are interpretable. Their results extend to sub-gaussian and SubWeibull cases also --- where they leverage the concentration inequalities in \citep{basu2015regularized}, we can simply substitute Proposition \eqref{DEandRE-SW},\eqref{DeviationSW-VAR} and \eqref{RE-SW-VAR}. The rest is identical.

\begin{figure}
%\centering
\caption{ \label{fig:figGL} \footnotesize Toy sparsity pattern (active groups shaded) induced by Own-Other Penalty.\\
\centerline{$p=3,d_A=4,d_B=2,s=9$}}
\includegraphics[scale=.6]{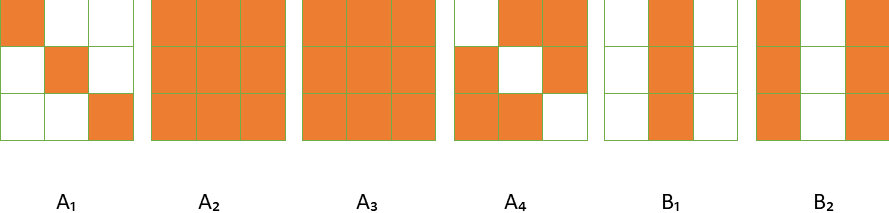}
\end{figure}

Formally, a $p$ dimensional centered VAR-X model $\{x_t\}$ with exogenous parts $\{z_t\}$, is given by
\begin{gather}
x_t=\sum_{i=1}^{d_{A}}A_i^\top x_{t-i}
+\sum_{j=1}^{d_{B}}B_j^\top z_{t-j} + \epsilon_t,
\end{gather}
with the exogenous process $\{z_t\}$ generated from a simple VAR(1) model $z_t= D^\top z_{t-1}+ \eta_t$, independent of the noise $\{\epsilon_t\}$. The autoregressive, endogenous lag is $d_A$, while the "distributed" exogenous lag is $d_B$. The endogenous and exogenous coefficients are respectively the $p$ by $p$ square matrices $A_i$ and $B_j$. Since we are primarily interested in sparsifying the parameters $A_i$, $B_j$, without loss, let the coefficient $D=0$. We can recast the VAR-X model as a VAR(1) model: $ y_t= F_{aug}^\top y_{t-1} + u_t$, as follows:
\begin{gather}\label{eq:VARX}
y_t:=
\begin{bmatrix}
x_t \\ 
x_{t-1} \\
\vdots \\
x_{t-d_A+1} \\
z_t \\
z_{t-1} \\
\vdots \\
z_{t-d_B}
\end{bmatrix}_{p(d_A+d_B)\times 1} ,
\quad
u_t:=
\begin{bmatrix}
\epsilon_t \\ 
0 \\
\vdots \\
0
\end{bmatrix}_{p(d_A +d_B)\times 1)}  ,
\quad 
F_{aug}^\top =\begin{bmatrix}
A_{aug}^\top & B_{aug}^\top\\
\end{bmatrix} 
\end{gather}
\begin{gather}
A_{aug}^\top :=
\begin{bmatrix}
A_1^\top & A_2^\top & \cdots & A_{d_A-1}^\top & A_{d_A}^\top \\ 
I_p & \bs{0} & \cdots & \bs{0}  & \bs{0} \\ 
\bs{0} & I_p &        & \bs{0}  & \bs{0} \\ 
\vdots &     & \ddots &  \vdots & \vdots\\ 
\bs{0} & \bs{0} & \cdots & I_p & \bs{0} \\
\bs{0} & \bs{0} & \cdots & \bs{0} & \bs{0} \\ 
\bs{0} & \bs{0} & \cdots & \bs{0} & \bs{0} \\ 
\bs{0} & \bs{0} &        & \bs{0} & \bs{0} \\ 
\vdots &     & \ddots &  \vdots & \vdots\\ 
\bs{0} & \bs{0} & \cdots & \bs{0} & \bs{0} 
\end{bmatrix}_{p(d_A+d_B)\times pd_A} 
\quad
B_{aug}^\top :=
\begin{bmatrix}
B_1^\top & B_2^\top & \cdots & B_{d_B-1}^\top & B_{d_B}^\top \\ 
\bs{0} & \bs{0} & \cdots & \bs{0}  & \bs{0} \\ 
\bs{0} & \bs{0}  &        & \bs{0}  & \bs{0} \\ 
\vdots &     & \ddots &  \vdots & \vdots\\ 
\bs{0} & \bs{0} & \cdots & \bs{0}  & \bs{0} \\
I_p & \bs{0} & \cdots & \bs{0} & \bs{0} \\ 
\bs{0} & I_p & \cdots & \bs{0} & \bs{0} \\ 
\bs{0} & \bs{0} &        & \bs{0} & \bs{0} \\ 
\vdots &     & \ddots &  \vdots & \vdots\\ 
\bs{0} & \bs{0} & \cdots & I_p & \bs{0} 
\end{bmatrix}_{p(d_A+d_B)\times pd_B} .
\end{gather}
Assuming the (augmented) noise $\epsilon_t$ is now a strictly stationary, mean-zero Subweibull process, we can carry out a similar consistency analysis as before, on the original process recast as $x_t=F^\top y_{t-1}+\epsilon_t$, where the parameter of interest is $F^\top=[A^\top:B^\top]$, by splitting it up into $p$ sub-processes running in parallel (see section 2). Interesting examples of the penalty $\mcR()$ are given in \citep[Table 1, page 7]{nicholson2017varx}.
All of them are variations of the (disjoint) group $\ell_1$ norm with different group structures with their own interpretations. To illustrate an example, consider the \textit{Own/Other} penalty
\begin{gather}
\mcR(F)=\mcR([A:B]):=\sqrt{p}\sum_{i=1}^{d_A} ||A_i^{\text{on}}||_F+\sqrt{p(p-1)}\sum_{i=1}^{d_A}||A_i^{\text{off}}||_F \\
+\sqrt{p}\sum_{j=1}^{d_B}\sum_{k=1}^p ||B_j^{(k)}||_F
\end{gather}
where $A_i^{\text{on}}$ and $A_i^{\text{off}}$ represent the vectors of diagonal and off-diagonal entries of $A_i$ respectively, and $B_j^{(k)}$ is the $k^{th}$ column of $B_j^\top$. Note that the penalty is weighted to avoid regularization favoring larger groups. A toy example is shown in Figure \eqref{fig:figGL} with the active (i.e. nonzero) elements shaded. This penalty is useful in many applications, where the diagonal entries of each $A_i$, which represent regression on a series' own lags, are  more likely to be nonzero than are off-diagonal entries, which represent lagged cross-dependence with other components. 

Clearly the parameter space is partitioned into $2d_A +pd_B$ disjoint groups with maximum group size $p(p-1)$, assuming only $s$ many groups are active, and the noise process $\epsilon_t$ has Subweibull norm $K$, Proposition \eqref{VAR-Subweibull-consistency} applies and gives the following estimates:
\begin{gather}
w(\mathbb{B}_\mathcal{R}(0,1))\le
[\sqrt{2d_A +pd_B}+2\sqrt{\log p(p-1)}]/\sqrt{p},\\
w^2(\mathcal{T}\cap\mathbb{B}_2)\le[(\sqrt{2\log(p(p-1)-s)} +\sqrt{2d_A +pd_B})^2+2d_A +pd_B](s/p),\quad\\
\Phi_\mathcal{R}(\mathcal{T})\le \sqrt{s},\quad
\bar{\Phi}\le 1.    
\end{gather}
Note that there is a factor of $p^{-1}$ due to the fact that the norm $\mcR()$ is weighted. This group $\ell_1$ LASSO results in the optimal consistency rate $\sqrt{\frac{s\log p}{n}}$, modulo a temporal dependence factor $\msfC(F_{aug})$.

\section{Numerical Experiments.}

\begin{figure}
%\centering
\caption{ \label{fig:figSW} \footnotesize The $\ell_1$-LASSO consistently estimates the Subweibull VAR transition matrix. \\
A smaller Subweibull tail index means heavier tails. }
\includegraphics[scale=.67]{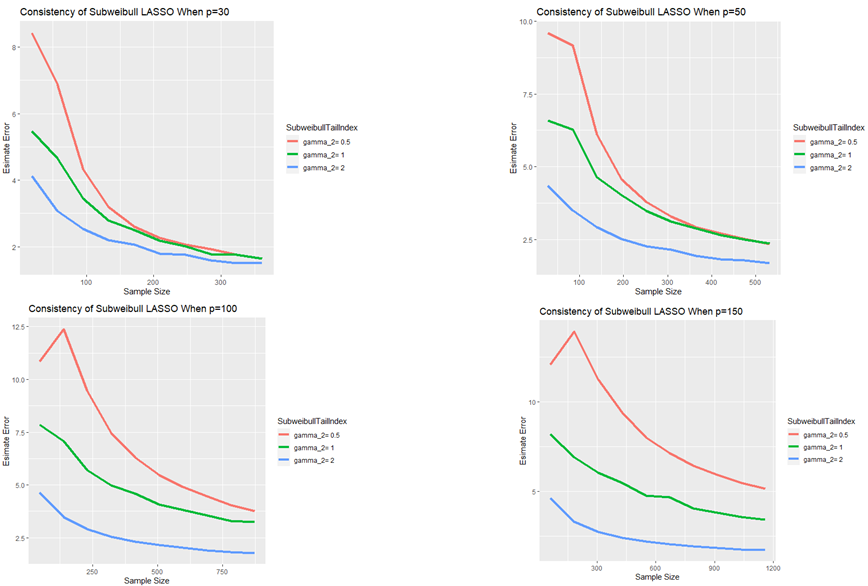}
\end{figure}

Following \cite{wong2020lasso}, we simulated $n$ samples from a $p$ dimensional VAR(1) model $z_t=B^\top z_{t-1}+\epsilon_t$, wherein the parameter matrix $B$ is element-wise sparse, having $s$ non-zero entries. The $p$ components of the noise $\epsilon_t$ are iid Weibull, with tail index $\gamma_2$. The transition matrix $B$ was simulated from an Erdos-Renyi graph with $p$ vertices and randomly choosing $s$ positions with nonzero entries and then sampling each non-zero entry in an iid fashion from a Uniform(0,1) distribution. Finally, we re-scaled the parameter $B$ to ensure its spectral radius is 0.5. We set $s=\sqrt{p}$ with $p\in\{30,50,100,150\}$ and select the following three regimes for $\gamma_2\in\{0.5,1,2\}$ that correspond to semi-exponential, exponential and sub-Gaussian tails, respectively. The sample size is set $n=m\times s\log p$ for $m\in\{1,3,5,7,9,13,15,17,19\}$. The estimated error $||\hat{B}-B||_F$, averaged over 30 replications, is plotted as a function of size $n$. Figure \eqref{fig:figSW} depicts the relationship between the Sub-Weibull tail index $\gamma_2$ and the estimated error. Note that a smaller $\gamma_2$ implies heavier tails, resulting in larger estimated error, as observed from the Figure.
\begin{figure}
\centering
\caption{\label{fig:figSW-tables}}
\includegraphics[scale=1.05]{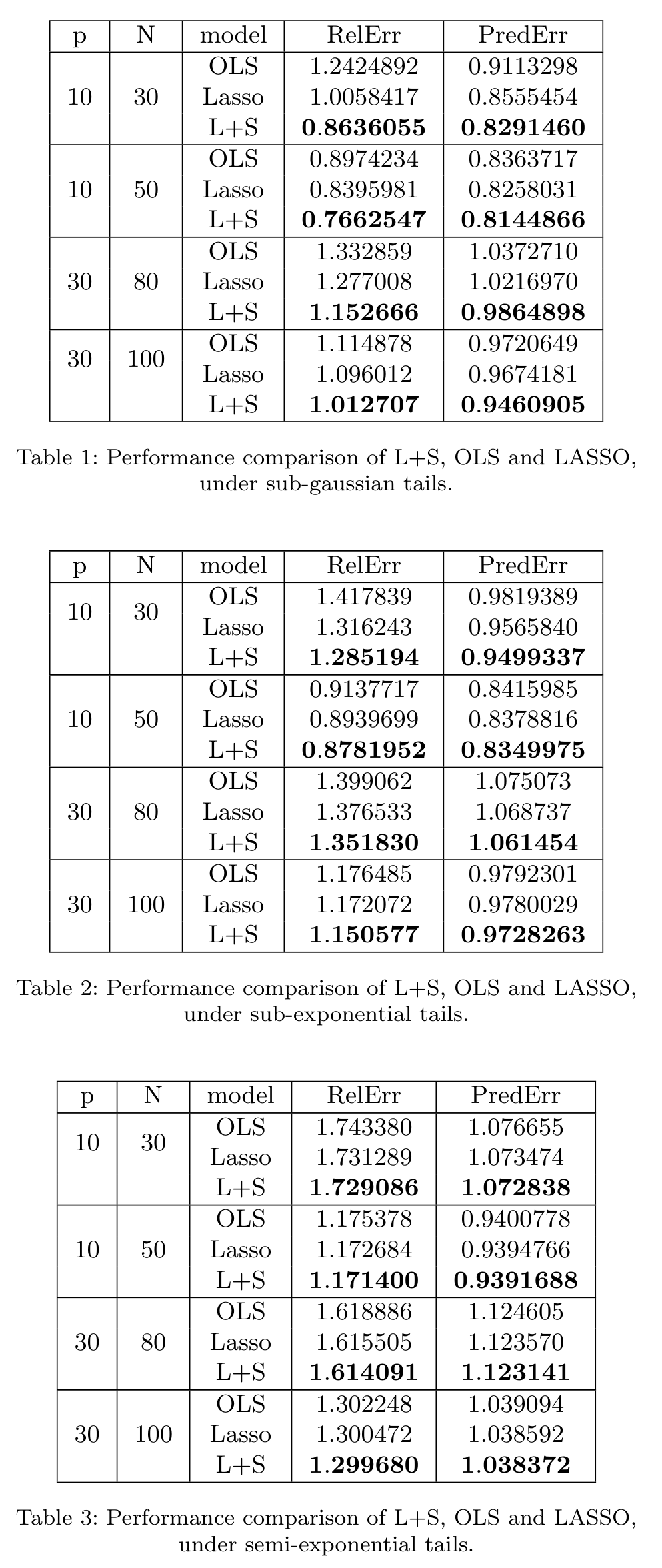}
\end{figure}  

Next, we estimate a "low rank + sparse" VAR(1) model, by using the following 3 methods --- (a) Ordinary Least Squares (OLS), (b) Sparse Lasso (LASSO) and (c) Low rank + sparse Lasso (LS). The dimensions and sample sizes are set as follows: (a) $p=10,n=30,50$ (small VAR), and (b) $p=30,n=80,100$ (medium VAR). The parameter transition matrix $B$ can be decomposed into a low-rank matrix $L$ (with fixed rank=3) and a sparse matrix $S$ with $2-4\%$ non-zero entries. We rescale the entries of $B$ to ensure stability of the process --- we set the spectral radius=0.7. 

We compute the relative estimated error (RelErr) and the out-of-sample prediction error (PredErr), under 3 regimes --- (a) $\gamma_2=2$ (sub-Gaussian), (b) $\gamma_2=1$ (sub-exponential) and (c) $\gamma_2=0.5$ (semi-exponential). We tabulate the results in Figure \eqref{fig:figSW-tables}. The number of out of samples used for calculating the prediction error is set to 10.

Under all 3 regimes of light-to-heavy tails, and both small and medium VAR, we find that the LS VAR estimates outperform the estimates using ordinary least-squares (OLS) and the sparse LASSO, since it produces the least prediction error, as expected. (Note: as the ratio $N/p$ increases, it is possible that the OLS may produce lower estimation errors than the LASSO in some cases, even if the OLS model is hard to interpret in high dimensions). Further, we note that the estimation errors produced from all 3 methods methods decrease with increasing sample sizes (given the dimension), as per theory. Finally the estimation errors are larger, when the tails are heavier, as expected from theory.

\section{Discussion}

In this paper, we consider the theoretical properties of penalized estimates in high-dimensional time series models when the data are generated from a (non-Gaussian) multivariate stationary process, under structured sparsity constraints. Specifically, we consider a sub-Weibull VAR$d$ model, and derived optimal consistency rates under general sparsity, and gave applications of the obtained results.

A interesting direction for future research is to carry out inference for the regression (transition matrix) parameter for a sub-Weibull model under structured sparsity constraints. For iid data and the LASSO penalty, a number of papers have obtained normal approximations for the parameter estimates after a debiasing steps - see, e.g, \cite{zhang2014confidence,van2014asymptotically,javanmard2014confidence,ning2017general}. The literature is rather sparse for temporally dependent data, with \cite{basu2019system} and  \citep{zheng2019testing} being two recent works on the topic. We believe that analogous results can be extended, with some care, to the case of Sub-Weibull VAR$d$ for debiased LASSO estimators. It would be of interest to extend such results to more general penalties $\mcR()$.

\begin{comment}
However, in some cases, we showed how to re-write the time series as Markov chains and get optimal consistency rates $\textit{and}$ optimal finite sample bounds. Still, there may be interesting examples of mixing processes that cannot be written in simple Markov form, or even if that were possible, it may be difficult to verify the ``drift'' and ``minorization'' conditions required to apply the tighter concentration bound \eqref{ConcMarkov} (e.g. multivariate GARCH models).
\end{comment} 
Finally, one potential future question which we briefly alluded to, is (non-asymptotic) inference (based on one-step or de-sparsified estimators) for high dimensional time series under general sparsity. 

%\appendix

\begin{appendices}

\section{Comparing the Concentration Bounds and Dependence Factors to Other Work:} \label{sec-Compare-Dependence}

\subsection{Comparison with Stability Factors in Basu et al.}

There are several notions of dependence for a stochastic process $\{x_t\}$. We only mention a few here. Basu et al \citep{basu2015regularized, basu2015network} quantify dependence in case of a strictly Gaussian time series by considering its spectral density. Formally, given a centered, stationary time series $\{x_t\}$ in $\mbR^p$, with autocovariance matrix function $\Sigma_X(h)= \cov(x_t,x_{t+h})$ (we write $\Sigma_X(0)$ as $\Sigma_X$ interchangeably, depending on context), its spectral density $f_X$ is a Hermitian matrix valued function given by:
\begin{gather}
f_X(\theta)=\frac{1}{2\pi} \sum_{l=-\infty}^{+\infty} \Sigma_X(l)e^{-il\theta} ,\quad \theta\in [-\pi,\pi].
\end{gather}
If $\{(x_t,y_t)\}$ is a jointly centered, stationary process with $\Sigma_{X,Y}(h)=\cov(x_t,y_{t+h})$, then their cross-dependence is given by the cross-spectral term
\begin{gather}
f_{X,Y}(\theta)=\frac{1}{2\pi} \sum_{l=-\infty}^{+\infty} \Sigma_{X,Y}(l)e^{-il\theta} ,\quad \theta\in[\pi,\pi] .  \end{gather}
\begin{definition}\label{stability}

Finally, assuming the existence of these densities, temporal dependence is quantified by the following stability factors:
\begin{gather}\label{stability_measures_defn}
\mcM(f_X):=\underset{\theta\in [-\pi,\pi]}{esssup} \quad\Lambda_{\max}(f_X (\theta))=\underset{\theta\in [-\pi,\pi]}{esssup}\quad ||f_X(\theta)||_2, \\
\mfm(f_X):=\underset{\theta\in [-\pi,\pi]}{essinf} \quad\Lambda_{\min}(f_X (\theta)) , \\
\mcM(f_{X,Y}):=\underset{\theta\in [-\pi,\pi]}{esssup} \quad\sqrt{\Lambda_{\max}f_{X,Y}^*(\theta) f_{X,Y}(\theta)} = \underset{\theta\in [-\pi,\pi]}{esssup} \quad ||f_{X,Y}(\theta)||_2 .
\end{gather}
\end{definition}

Also, for $1\le k\le p$, the factors $\mcM(f_X,k)= \underset{J\subseteq \{1,...,p\}, card(J)\le k}{max} \mcM(f_{X_J})$ and similarly $\mfm(f_X, k)$, may be defined. A large class of time series (VAR models in particular) satisfies the assumption $0<\mfm(f_X)\le\mcM(f_X)<\infty$. The term $\mcM(f_X)$ measure the "spikes" of the spectrum which indicate greater dependence or "memory" in the process, while the term $\mfm(f_X)$ measures the cross sectional dependence between the $p$ component processes. To estimate stability factors, consider a linear filter of an innovation process $\{\eta_t\}\sim\text{IID}(0,\Sigma_\eta)$, given by $x_t=\mcA(\mathfrak{B})\eta_t$ where $\mathfrak{B}$ denotes the backshift operator, the matrix power series $\mcA(z)=\sum_{i\ge 0}A_i z^i$ defined on the complex plane $z\in\mbC$, satisfy $\det(\mcA(z))\neq 0$ on the unit disk $\{z\in\mbC:|z|\le 1\}$  (stationarity), and $\sum_{i\ge 0}||A_i||_2<\infty$ (stability). Its spectral density is given by 
\begin{gather}
f_X(\theta)=\mcA(e^{-i\theta})\frac{\Sigma_\eta}{2\pi}
\mcA^*(e^{-i\theta}) ,\quad \theta\in [-\pi,\pi] ,
\end{gather}
which leads to the following estimates
\begin{gather}\label{meas_stability_est}
\frac{\Lambda_{\min}(\Sigma_\eta)\mu_{min}(\mcA)}{2\pi}\le \mfm(f_X)\le\mcM(f_X)\le \frac{\Lambda_{\max}(\Sigma_\eta)\mu_{max}(\mcA)}{2\pi} ,
\end{gather}
where
\begin{gather}
\mu_{max}(\mcA):=\underset{|z|=1}{\max} \quad \Lambda_{\max}(\mcA(z)^*\mcA(z))
=\underset{|z|=1}{max}||\mcA(z)||_2^2 ,\\
\mu_{min}(\mcA):=\underset{|z|=1}{\min}\quad \Lambda_{\min}(\mcA(z)^*\mcA(z))
=\underset{|z|=1}{\min}||\mcA(z)||_2^{-2} .
\end{gather}
This can easily be extended to a stable rational filter of the form $\mcA/\mcB$ in which case, we will get bounds matching (2.6) in \citep{basu2015regularized}. The advantage of using these temporal dependence factors is that (i) they are easy to interpret and (ii) we can state tight concentration inequalities for dependent Gaussian processes, in terms of these factors, as stated below.

\begin{proposition}\label{ConcBasu}
{\citep*[Proposition 2.4]{basu2015regularized}.}

Let $\{x_t\}$, $\{y_t\}$ be centered, stable Gaussian processes with $x_t$ and $y_t$ being independent for each $t$. Let $v\in\mbR^p$ be a unit vector, i.e. the Euclidean norm $||v||\le 1$. Then
\begin{gather}
\mbP[|\sum_{i=1}^n v^\top x_i y_i^\top v| >2\pi \parenc{\mcM(f_X)+\mcM(f_Y)+\mcM(f_{X,Y})}nt] \le 6\exp[-cn\min\parenc{t,t^2}] , \\
\text{and,} \quad
\mbP[|\sum_{i=1}^n v^\top[x_i x_i^\top -\mbE(x_1 x_1^\top)]v| >2\pi\mcM(f_X)nt] \le 2\exp[-cn\min\parenc{t,t^2}] ,
\end{gather}
for $t>0$ and an absolute constant $c>0$. This is just the Hanson-Wright inequality \citep{vershynin2018high} modulo temporal dependence factors $\mcM(f_X)+\mcM(f_Y)+\mcM(f_{X,Y})$ and $\mcM(f_X)$.

\end{proposition}

As noted, our concentration bound Proposition \eqref{DEandRE-SW} directly compares with Proposition \eqref{ConcBasu} for the strictly Gaussian case. Finally, this dependence based on the spectral density $f_X$, can be extended to linear time series with Subweibull tails. Specifically, we can connect the measure of dependence defined in section \eqref{SubweibullVAR}, and the stability factors, as follows.

\begin{proposition}

We have that
\begin{gather}
\Lambda_{\min}(\Sigma_\eta) ||A_0||_2^2 \parens{\msfC(\mcA)}^{-2}\le2\pi\mfm(f_X)
\le 2\pi\mcM(f_X)
\le 2\Lambda_{\max}(\Sigma_\eta)\msfC(\mcA).
\end{gather}
Often, $A_0=I_p$ by convention, so $||A_0||_2 =1$.

\end{proposition}

\begin{proof}

First, $\forall\theta\in[-\pi,\pi]$,

\begin{align}
2\pi f_X(\theta)
=\Sigma_X(0)+ \sum_{l=1}^\infty [\Sigma_X(l)e^{-il\theta} +\Sigma_X(-l)e^{il\theta} ]\\
=\Sigma_X(0)+ \sum_{l=1}^\infty [\Sigma_X(l)e^{-il\theta} +\Sigma_X(l)^\top e^{il\theta} ]\\
\text{so,}\quad
2\pi\mcM(f_X)\le 2\sum_{l\ge 0}||\Sigma_X(l)||_2
\le 2\Lambda_{\max}(\Sigma_\eta) \sum_{i,j\ge 0} ||A_{j+i}||_2||A_j||_2 \\
= 2\Lambda_{\max}(\Sigma_\eta)\msfC(\mcA).
\end{align}

To find a lower bound, we note that, from section \eqref{temp-dep},
\begin{gather}
2\pi\mfm(f_X)\ge\Lambda_{\min}(\Sigma_\eta) \mu_{min}(\mcA),\\
\text{where,}\quad\mu_{\min}(\mcA) =\underset{|z|=1}{\min}||\mcA(z)||_2^{-2} 
\ge\parens{\sum_{i\ge 0}||A_i||_2}^{-2}
\ge ||A_0||_2^2 \parens{\msfC(\mcA)}^{-2} ,\\
\text{so,}\quad 2\pi\mfm(f_X)\ge \Lambda_{\min}(\Sigma_\eta) ||A_0||_2^2 \parens{\msfC(\mcA)}^{-2} .
\end{gather}
\end{proof}

\begin{remark}
Basu et al specifically assume Gaussianity to quantify dependence and derive their concentration bounds. It is expected that their stability factors would extend to heavier tails but this is not obvious, since the trick they exploit in the Gaussian case --- rotating a Gaussian vector to achieve independence--- fails even for the \textit{sub}-Gaussian case, e.g. when the data is discrete or bounded. Our dependence factor directly compares with theirs, with the added advantage that ours is also applicable in a SubWeibull setting. Nonetheless, for completeness, we derive consistency rates for Gaussian stochastic regressions under a general sparsity structure in Appendix \eqref{sec-Estimation-proc}. %\textcolor{red}{WHICH APPENDIX??}
\end{remark}

\subsection{Comparison with Wong et al: Mixing.}

Wong \citep{wong2016lasso,wong2020lasso} considers more classical notions of dependence, specifically that of $\textit{mixing}$. For a probability space $(\Omega,\mcF,\mbP)$, given two sub $\sigma$ fields $\mcA,\mcB\subseteq\mcF$, the following notions of mixing \citep{bradley2005basic} are defined:
\begin{gather}
\alpha(\mcA,\mcB):=\sup\{|\mbP(A\cap B) -\mbP(A)\mbP(B)|,\quad A\subseteq\mcA, B\subseteq\mcB\} ,\\
\rho(\mcA,\mcB):=\sup\{\text{Corr}(f,g),\quad f\in\mcL_{real}^2(\mcA),g\in\mcL_{real}^2(\mcB) \} ,\\
\beta(\mcA,\mcB):=\frac{1}{2}\sup[\sum_{i=1}^I \sum_{j=1}^J |\mbP(S_i\cap T_j)-\mbP(S_i)\mbP(T_j)| ,\\
\Omega=\sqcup_{i=1}^I A_i=\sqcup_{j=1}^J B_j,\quad A_i\in\mcA,B_j\in\mcB ]
\end{gather}
where $\sqcup$ denotes a disjoint union, and $\mcL_{real}^2(\mcA)$ denotes the space of square integrable, $\mcA$ measurable, real valued random variables. We have the following estimates: $2\alpha(\mcA,\mcB)\le\beta(\mcA,\mcB)$,
$4\alpha(\mcA,\mcB)\le\rho(\mcA,\mcB)$  (p109,\citep{bradley2005basic}). If $\mcA=\sigma(X)$, $\mcB=\sigma(Y)$ for two random elements $X$,$Y$, then we write $\alpha(X,Y)$ instead of $\alpha(\mcA,\mcB)$, etc. For a (strictly) stationary random sequence $\{(x_t,y_t)\}$, denoting $x_{-\infty:t}=\sigma\{x_j,j\le t\}$, etc. we define, for $l\in\mbZ$ 
\begin{gather}
\alpha_X(l)=\alpha(x_{-\infty:t},x_{t+l:\infty}),\quad
\alpha_{X,Y}(l)=\alpha(x_{-\infty:t},y_{t+l:\infty}),\quad 
\text{etc.}
\end{gather}
The usual practice is to define the mixing coefficients for lags $l\ge 1$. We will also use the following fact repeatedly (especially in the non-gaussian case): since $\mcA'\subseteq\mcA$, $\mcB'\subseteq\mcB$ implies $\alpha(\mcA',\mcB')\le\alpha(\mcA,\mcB)$ (same for $\beta$ or $\rho$ mixing), for a mixing sequence $\{x_t\}$, and a $\textit{finite}$ lag $d$, the process $\{y_t\}$ defined by any measurable function $y_t:=g(x_t,x_{t-1},\dots,x_{t-d})$, is also mixing (in the same sense as $\{x_t\}$), with mixing coefficients bounded above by those of the original sequence $\{x_t\}$ (see e.g. \citep*[Theorem 14.1]{davidson1994stochastic}).

\begin{definition}

The process $\{x_t\}$ is said to be $\alpha$-mixing (or $\textit{strong}$ mixing) if $\alpha_X(l)\rightarrow 0$, $\beta$-mixing (or $\textit{regular}$) if $\beta_X(l)\rightarrow 0$ and $\rho$-mixing (or $\textit{uniformly}$ mixing) if $\rho_X(l)\rightarrow 0$, as $l\rightarrow\infty$. Moreover a $\beta$-mixing process $\{x_t\}$ is $\textit{sub}$ geometric if there are positive constants $\beta_{mix}<1$ and $\gamma>0$ such that
\begin{gather}
\beta(n)\le 2\exp{(-n^\gamma\log\beta_{mix}^{-1})} \quad \text{for all $n\in\mbN$ .}
\end{gather}
We call this $\gamma$ its geometric index, and $\beta_{mix}$ the mixing rate (note that $\beta_{mix}$ may depend on the true regression parameter $\beta^*$ if $\{x_t\}$ is endogenous). In applications, usually $\gamma\le 1$. If $\gamma=1$, $\{x_t\}$ is said to be $\textit{geometrically}$ $\beta$-mixing. 

\end{definition}

Several more mixing notions exist in literature. For a connection between properties of spectral density and mixing, see \citep[sections 6,7]{bradley2005basic}. In particular, for a stationary, centered Gaussian sequence $\{(x_t,y_t)\}$, we can connect the stability factors to the mixing coefficients as follows:
\begin{proposition}
\begin{gather}
\mcM(f_X)\le||\Sigma_x||_2\sum_{l\in\mbZ}\alpha_X(l) ,\quad \mcM(f_{X,Y})\le\sqrt{||\Sigma_x||||\Sigma_y||}\sum_{l\in\mbZ} \alpha_{X,Y}(l) .   
\end{gather}
\end{proposition}

This follows from \citep[page 32]{wong2020lasso} and the equivalence of $\alpha$ and $\rho$ mixing for Gaussian processes \citep[Theorem 2]{kolmogorov1960strong}. In particular, if the mixing coefficients $\{\alpha_X(l)\}$ (or equivalently $\{\rho_X(l)\}$) are summable, then $\mcM(f_X)<\infty$. This is satisfied for many Gaussian processes. Further, denoting the mixing coefficients of the joint process $\{(x_t,y_t)\}$ by $\alpha_{(X,Y)}$ (not to be confused with $\alpha_{X,Y}$), it follows that $\alpha_X,\alpha_Y, \alpha_{X,Y}\le\alpha_{(X,Y)}$. \citep{wong2020lasso}, among others, prefers to work with the stronger assumption that the joint process $\{(x_t,y_t)\}$ is mixing, hence relies on coefficients $\alpha_{(X,Y)}$. Note that, in the non Gaussian case, $\alpha$ and $\rho$ mixing are usually not equivalent and $\rho$ mixing is often too restrictive for many time series. In those cases, $\beta$-mixing, though stronger than $\alpha$-mixing, is an acceptable compromise. It would be interesting to connect the stability factor $\mcM(f_X)$ with some decay condition on the mixing coefficients, in those cases. For multivariate processes (which is of interest in our case), theoretical results connecting spectral properties with the mixing properties of stationary processes are sparse in the literature
\citep{cheng1993mixing}.  

There exist concentration inequalities in terms of mixing coefficients in literature. For the Gaussian case, we can write a concentration bound similar to \citep*[Proposition 2.4]{basu2015regularized} by replacing the stability factors by the mixing coefficients (this leads to \citep*[Propositions 2 and 3]{wong2020lasso}). Going beyond Gaussianity poses significant problems to the mixing framework, since tight concentration bounds are not easily available. Using these bounds, in turn, leads to suboptimal consistency rates or suboptimal minimum sample size required to achieve that rate (with high probability). The key concentration inequality used in this case, holds for Subweibull random variables under $\beta$-mixing, and is given by \citep[lemma 13]{wong2020lasso}, which itself is a modification of a quite general result in \citep{merlevede2011bernstein} (such results are few and far in between).

\begin{proposition} \label{ConcWong}
{Lemma 13, \citep{wong2020lasso}.}

Let $\{x_t\}$ be a centered, strictly stationary, one dimensional, Subweibull($\gamma_2$/2) process with Subweibull norm bounded by $K$. Further, assume that it is $\textit{sub}$-geometrically $\beta$-mixing with mixing rate $\beta_{mix}$ and sub-geometric index $\gamma_1$. Let $n\ge 4$ and $t>1/n$. Also, let $1/\gamma_1+2/\gamma_2=1/\gamma$ and suppose $\gamma<1$. Then,
\begin{gather}
\mbP\parenf{|\sum_{t=1}^n x_t|>Knt}\le n\exp\parens{-\frac{1}{\msfC^2} \min\parenc{(nt)^\gamma,nt^2}} , 
\end{gather}
wherein $\msfC>0$ is a constant that depends on $\beta_{mix}$, $\gamma_1$ and $\gamma_2$. 

\end{proposition}

\begin{remark}
In Proposition \eqref{ConcWong}, we specifically assumed the process $\{x_t\}$ has a SubWeibull tail index $\gamma_2/2$, and not $\gamma_2$, in order to make a direct comparison to Proposition \eqref{DEandRE-SW}. In the latter case, the Subweibull terms are essentially squared, which leads to their tail index being halved \citep[Lemma 6]{wong2020lasso}.Hence the term $2/\gamma_2$ appears in the right hand side of the concentration inequalities in \eqref{DEandRE-SW}. By comparison, in the right hand side of Proposition \eqref{ConcWong}, we find the term $\gamma$ instead of $\gamma_2/2$. Since this $\gamma<\gamma_2/2$, this concentration bound is suboptimal. Another drawback of Proposition \eqref{ConcWong} is that the dependence factor $\msfC^{-2}$ is intractable compared to ours, even for the special case of SubWeibull VAR.   
\end{remark}

\begin{remark}
As noted, this inequality does not hold for all $t>0$. However, the precondition $t>1/n$ is extremely mild and boils down to the sample size $n$ being bigger than a finite constant (that does not scale with dimension) in applications. Also, a factor of $n$ still remains on the right hand side, which is sub-optimal. However, this is not a big price to pay. In fact, the exponential term usually dominates. For completeness, we derive consistency rates for SubWeibull stochastic regressions in Appendix \eqref{Stoch-reg-SW}. 
\end{remark}

\section{Further Concepts Related to the Penalty Norm.}
\label{sec-Futher-Concepts-Penalty}

Here, we present some general concepts related to sparsity inducing norms, since many interesting norms in practice are $\textit{not}$ decomposable (unlike the $\ell_1$ norm) and require additional technical machinery.

\begin{definition}
A norm $\mcR$ is called $\textit{decomposable}$, if for any $S\subseteq\{1,...,p\}$, for all $v\in\mbR^p$, we have $\mcR(v)=\mcR(v_S)+\mcR(v_{S^c})$.
\end{definition}

\begin{definition}
Following \citep{chandrasekaran2012convex}, let $\mcA$ be a collection of $\textit{atoms}$ that is a compact subset of $\mbR^p$. Suppose no element $a\in\mcA$ lies in the convex hull of the other elements $conv(A-\{a\})$, i.e., the elements of $\mcA$ are the extreme points of $conv(A)$. Let $||x||_\mcA$ denote the gauge of $\mcA$, i.e. $||x||_\mcA= \inf\{t>0 :x\in t\cdot conv(\mcA)\}$.The gauge is always a convex, extended-real valued function for any set $\mcA$. By convention this function evaluates to $+\infty$ if $x$ does not lie in the affine hull of $conv(\mcA)$. We will assume without loss, that the centroid of $conv(\mcA)$ is at the origin, as this can be achieved by appropriate re-centering. With this assumption the gauge function may be recast as:
\begin{gather} \label{atomic-norm}
||x||_\mcA=\inf\{\sum_{a\in\mcA} c_a: x=\sum_{a\in\mcA}c_a a,c_a\ge 0, \forall a\in\mcA \} ,   
\end{gather}
with the sum being replaced by an integral when $A$ is uncountable. If $A$ is centrally symmetric about the origin (i.e., $a\in\mcA$ if and only if $-a\in\mcA$) we have that $||\cdot||_\mcA$ is a norm, which we call
the $\textit{atomic}$ norm induced by $\mcA$. Henceforth, we assume that $||\cdot||_\mcA$ is a norm. These norms are extremely useful because they need not be decomposable, but still cover a rich class of norms (see \citep{chandrasekaran2012convex,bhaskar2013atomic}).

\end{definition}

\begin{remark}

We can also induce sparsity via non-convex penalties like the SCAD and MCP \citep{fan2001variable,zhang2010nearly}. It is known that in certain high-dimensional regimes, the estimation error of nonconvex penalties like SCAD, MCP scales roughly in the same order as LASSO. Of course these results are established with iid data. However, they can easily be extended to our setup as well, since only the algebraic arguments need to be modified --- the probabilistic arguments remain the same.

\end{remark}

\section{Stochastic Regression under Gaussianity.}
\label{sec-Estimation-proc}
In this section, we consider the stochastic regression model \eqref{stoch-reg-matrix-form} under Gaussianity  and the previously defined estimation procedures ---the Dantzig selector and the LASSO--- subject to generic regularizers $\mcR(\cdot)$. First, we derive consistency results for the Dantzig selector --- the results for LASSO are analogous. 

\subsection{The Dantzig selector}

The Dantzig selector defined in \eqref{dantzig-selector-def} has been extensively studied in literature for $\mcR()=||\cdot||_1$ and then extended at substantial depth \citep*{chandrasekaran2012convex, chen2015structured,chen2016structured}, to the case of general regularizers, for iid data. To ensure consistency, a first and a second order conditions need to hold with high probability \citep*{wainwright2019high}. Specifically, the quantity $\mcR^*\parenf{X^\top\mcE/n}$ needs to concentrate around zero (\textit{deviation condition}), and also the minimum eigenvalue of the sample Gram matrix $X^\top X/n$ restricted to a small set, needs to be bounded uniformly away from zero (\textit{Restricted Eigenvalue} (RE) condition). These conditions are stated precisely below:

\begin{proposition}{Deviation condition:}\label{dev-stochreg}

There is a constant $c_0>0$ such that, for $n\ge c_0 w^2[\mbB_\mcR(0,1)]$, 
\begin{gather}\label{deviation-bound}
\mbP\left[ \mcR^*\parenf{\frac{X^\top\mcE}{n}}\ge 4\pi[\bar{\Phi}_\mcR^2 \mcM(f_X,1) + \mcM(f_\epsilon) + \bar{\Phi}_\mcR \mcM(f_{X,\epsilon})]\sqrt{\frac{c_0 w^2(\mbB_\mcR (0,1))}{n}} \right]  \\
\le 6\exp[-w^2(\mbB_\mcR(0,1))],
\end{gather} 
where $w[\mbB_\mcR(0,1)]$ is the Gaussian width of the unit norm ball $\mbB_\mcR(0,1)$. The terms $\mcM(f_X,1)$, $\mcM(f_\epsilon)$ and $\mcM(f_{X,\epsilon})$ quantify the temporal dependence due to a one-dimensional projection of the covariate process $\{x_t\}$, the noise process $\{\epsilon_t\}$, and the cross-dependence between the two, respectively ---see Appendix \eqref{normandtempdep}.

\end{proposition}

\begin{proposition}{RE condition.}\label{RE-stochreg}

Suppose $\mfm(f_X)>0$ and $c_0'>0$ is a constant. Let the tangent cone at $\beta^*$ be denoted by $\mcT=\mcT_\mcR(\beta^*)= \mathrm{cone}\parenf{\{v:\mcR(v+\beta^*)\le\mcR(\beta^*)\}}$. Then, a sample size of $n\ge 16c_0' \parenc{\frac{\mcM(f_X)}{\mfm(f_X)}}^2 w^2(\mcT\cap\mbB_2)$ is sufficient to guarantee that
\begin{gather} \label{RE-ineq}
\underset{v\in\mcT\cap\mbB_2}{inf}\frac{v^\top X^\top X v}{n} \ge\alpha_{RE} ,
\end{gather}
with probability at least $1-2\exp[-w^2(\mcT\cap\mbB_2)]$, where the restricted eigenvalue is $\alpha_{RE}=\pi\mfm(f_X)$. Further, $\mfm(f_X)$ may be replace by $\Lambda_{\min}(\Sigma_X)$/2, if required.

\end{proposition}

\begin{remark}
The deviation condition implies that, under the scaling $w^2(\mbB_\mcR (0,1))/n\rightarrow\infty$, the deviation term $\mcR^*(X^\top\mcE/n)$ concentrates around zero. The RE condition, on the other hand requires the sample size $n$ to scale with the size of the spherical cap of the tangent cone $\mcT$. These results are standard by now in case of iid data,
the only difference in our setting being the price paid for temporal dependence. Specifically, we consider the dependence factors appearing in the deviation and RE conditions. Ignoring the reverse compatibility factor $\bar{\Phi}_\mcR$ (it is usually bounded by an absolute constant), the dependence factor from the deviation condition is $\mcM(f_X,1) +\mcM(f_\epsilon)+\mcM(f_{X,\epsilon})$. This dependence factor appears in Proposition 3.2, \citep*{basu2015regularized}. It also matches the dependence factor in Proposition 2, \citep*{wong2020lasso}. To see this, note that, under the assumption the Gaussian process $\{(x_t,y_t)\}$ (equivalently  $\{(x_t,\epsilon_t)\}$ by linearity) is mixing with coefficients $\{\alpha(l)\}$ that sum to $\alpha$,
\begin{align}
\mcM(f_X,1)+\mcM(f_\epsilon)+\mcM(f_{X,\epsilon}) 
&\le\alpha\parenc{||\Sigma_X||_2+\sigma_\epsilon^2 + \sqrt{||\Sigma_X||_2 \sigma_\epsilon^2} } \\
&\le\frac{3\alpha}{2}\parenf{||\Sigma_X||_2+\sigma_\epsilon^2} 
\quad\text{(A.M.-G.M. inequality)} \\
&\le\frac{3\alpha}{2}\parenf{||\Sigma_X||_2+2\sigma_y^2+
2\beta^{*\top}\Sigma_X\beta^* } \\
&\le 3\alpha\parenf{\sigma_y^2+(1+||\beta^*||)||\Sigma_X||_2} .
\end{align}

However, the dependence factor from the RE condition is $\parenc{\mcM(f_X)/\mfm(f_X)}^2$, which is worse than that appearing in Proposition 3.1, \citep*{basu2015regularized}, since $\mcM(f_X, 2k)$ can be much smaller than  $\mcM(f_X)$, if $k$ (which denotes the sparsity of $\beta^*$) is much smaller than $p$. However, this is due to the decomposability of the $\ell_1$ norm, which gives rise to a cone set that can be $\textit{directly}$ approximated by (contained within) a $k$-sparse set (lemma F.1), whose Gaussian width is easy to compute. However, this is not the case for non-decomposable norms. In case of $\textit{endogenous}$ processes, such as VAR, the spectrum $f_X$ depends on the sparse $\beta^*$; then, we can estimate the stability factors in terms of the low dimensional parameter $\beta^*$, and thus control them. Note that, in case of $\alpha$-mixing, we have (using the above notation), $\mcM(f_X)\le \alpha||\Sigma_X||_2 =\alpha\Lambda_{\max}(\Sigma_X)$, and, as remarked in Proposition \eqref{RE-stochreg}, replace $\mfm(f_X)$ by $\Lambda_{\min}(f_X)/2$. Hence, we can replace the factor $\parenc{\mcM(f_X)/\mfm(f_X)}^2$ by $\parenc{\alpha\Lambda_{\max}(\Sigma_X)/ \Lambda_{\min}(\Sigma_X)}^2$, which matches the factor $``\eta^{-2}"$in Proposition 3, \citep*{wong2020lasso}.

\end{remark}

\begin{proposition}{Theoretical consistency of the Dantzig selector.}\label{det-errbd-dant-stochreg}

Suppose that the tuning parameter $\lambda_n$ and the sample size $n$ satisfy the following bounds:
\begin{gather}
\lambda_n = 4\pi[\bar{\Phi}_\mcR^2 \mcM(f_X,1) + \mcM(f_\epsilon) + \bar{\Phi}_\mcR \mcM(f_{X,\epsilon})\sqrt{\frac{c_0 w^2(\mbB_\mcR (0,1))}{n}} , \\
n\ge\max\parenc{c_0 w^2(\mbB_\mcR (0,1)), 16c_0' \parenc{\frac{\mcM(f_X)}{\mfm(f_X)}}^2 w^2(\mcT\cap\mbB_2)} \\
\succsim\parenc{\frac{\mcM(f_X)}{\mfm(f_X)}}^2\max\parenc{w^2(\mbB_\mcR(0,1)),w^2(\mcT\cap\mbB_2)}.
\end{gather}
Then, the penalized estimate $\hat{\beta}$ satisfies:
\begin{gather}
||\hat{\beta}-\beta^*||\leq \frac{2\lambda_n\Phi_\mcR(\mcT)}{\alpha_{RE}}, \quad
\mcR(v)\le\frac{2\lambda_n\Phi_\mcR^2(\mcT)}{\alpha_{RE}},
\quad\text{(Estimation error)}, \\
\frac{v^\top X^\top Xv}{n}\le \frac{4\lambda_n^2\Phi_\mcR(\mcT)}{\alpha_{RE}} 
\quad\text{(Prediction error)}.
\end{gather} 
with probability at least $1-6\exp[-w^2(\mbB_\mcR(0,1))]-2\exp[-w^2(\mcT\cap\mbB_2)]$. Here $c_0$ and $c_0'$ are constants that appear in the deviation and RE conditions respectively.
 
\end{proposition}

\subsection{The LASSO}
\label{LASSO}

The LASSO type estimator is the solution to the following constrained optimization problem:
\begin{gather}
\hat{\beta} = \underset{\beta\in\mbR^p}{
\mathrm{\argmin}} \quad \frac{1}{n} ||Y - X\beta||^2 + \lambda_n \mcR(\beta) .
\end{gather}
Here the penalty is given by a generic norm $\mcR()$ which reflects our beliefs about the underlying structure about the parameter $\beta^*$ in the model (\eqref{stoch-reg-matrix-form}). Also, $\lambda_n$ is a tuning parameter. This estimator has been shown to be consistent in the context of Gaussian stochastic regression, using $\ell_1$ norm \citep*{basu2015regularized}. We extend the result from $\ell_1$ to a general penalty $\mcR$. As we will see, the LASSO and Dantzig are virtually equivalent, with the same consistency rates, for any norm $\mcR$. The only difference is the cone of anti-concentration $\mcT$, given in Proposition \eqref{RE-stochreg}, will be replaced by a larger cone $\msC$ given by
\begin{gather}
\msC=\msC(\beta^*) :=\mathrm{cone}\{u:\mcR(u)/2 +\mcR(\beta^*) -\mcR(\beta^*+u)\geq 0\} .    
\end{gather}
In case $\mcR$ is decomposable, $\mcT\subseteq \msC\subseteq\mfC(J,3)$, where $\mfC(J,3)=\{v:\mcR(v_{J^c}) \le 3\mcR(v_J)\}$ is the standard cone introduced in \citep*{bickel2009simultaneous}. We may prove Proposition \eqref{RE-stochreg} with $\mcT$ replaced by $\msC$. Thus, we first give deterministic recovery bounds for the LASSO.

\begin{proposition}\label{det-errbd-LASSO-stochreg}

Assume the following:
\begin{gather}\label{deviation-algebra}
\lambda_n \geq 4 \mcR^*\parenf{\frac{X^\top \mcE}{n}} \quad \text{(Deviation Condition)},
\end{gather}
and that, there exists and $\alpha_{RE}>0$ such that 
\begin{gather}\label{RE-algebra}
\underset{v\in\msC\cap\mbB_2}{inf} \frac{1}{n} ||Xv||^2 \geq\alpha_{RE} \quad \text{(RE condition)}  
\end{gather}
Then, the penalized estimate $\hat{\beta}$ satisfies:
\begin{gather}
||\hat{\beta}-\beta^*|| \leq \frac{3\lambda_n\Phi_\mcR(\msC)}{2 \alpha_{RE}} , \\
\frac{||X(\hat{\beta}-\beta^*)||^2}{n}\le \frac{9\lambda_n^2\Phi_\mcR^2(\msC)}{2\alpha_{RE}} .
\end{gather} 

\end{proposition}

Of course, the deviation and RE conditions can be shown to hold with high probability as in case of the Dantzig selector. The tuning parameter and sample size scales similarly.

\begin{remark}

Of course, all our results carry through for a multi-response regression $y_t^\top=x_t^\top B +\epsilon_t^\top$ where $\{x_t\}$ and $\{y_t\}$ are now random processes in $\mbR^p$ and $\mbR^q$ respectively, and $B$ is a $p\times q$ matrix of regression parameters. We simply vectorize the process as follows:
\begin{gather}
\underbrace{
\begin{bmatrix}
y_1^\top \\ 
y_2^\top \\
\vdots \\
y_n^\top
\end{bmatrix}}_{Y} 
=\underbrace{
\begin{bmatrix}
x_1^\top \\ 
x_2^\top \\
\vdots \\
x_n^\top
\end{bmatrix}}_{X} \cdot B +
\underbrace{
\begin{bmatrix}
\epsilon_1^\top \\ 
\epsilon_2^\top \\
\vdots \\
\epsilon_n^\top
\end{bmatrix}}_{E}
\end{gather}
So that the matrix form of the process can be written in vectorized form: $vec(Y)=(I_q\otimes X)vec(B) +vec(E)$.

\end{remark}

\section{Stochastic regression with Subweibull tails and mixing.}\label{Stoch-reg-SW}

The concentration \eqref{ConcWong} compares favorably with the Gaussian case \citep[Proposition 2.4]{basu2015regularized}, \citep[lemma 11]{wong2020lasso}), but for the factor $n$ in the right hand side. Also, as it will be clear from the proofs that utilize this result, the precondition $t\ge 1/n$ implies $n$ is bigger than a finite constant (as opposed to scaling with dimension $p$, or quantities related to $p$). Hence, this will not be mentioned explicitly in statements giving finite sample bounds in terms of dimension $p$, to avoid clutter. For a more explicit derivation, see e.g \citep*[Proposition 8]{wong2020lasso}.

Suppose now we have the regression model in Section \eqref{sec-model-and-background} but instead of the processes
$\{x_t\}$, $\{\epsilon_t\}$ (and hence the response $\{y_t\}$) being Gaussian, we make the following probabilistic assumptions about the joint process $\{(x_t,\epsilon_t)\}$:

\begin{itemize}
    \item The process $\{(x_t,\epsilon_t)\}$ is centered and strictly stationary. Also, $\mbE[\epsilon_t|x_t]=0$ for each $t$.
    \item The process $\{(x_t,\epsilon_t)\}$ is Subweibull($\gamma_2$) with $||(x,\epsilon)||_{\psi_{\gamma_2}}\le K$.
    \item The process $\{(x_t,\epsilon_t)\}$ is (sub)geometrically  $\beta$-mixing with rate $\beta_{mix}$, and exponent $\gamma_1$.
    \item We  assume $\gamma<1$, where
    \begin{gather}
        \gamma:=\parens{1/\gamma_1+2/\gamma_2}^{-1}.
    \end{gather}
\end{itemize}

\begin{remark}
By\citep[Fact 1]{wong2020lasso} and linearity, assuming $\{(x_t,\epsilon_t)\}$ is mixing is equivalent to assuming $\{(x_t,y_t)\}$ is mixing.
\end{remark}

\begin{remark}
In applications, usually $\gamma_2\le 2$. In that case, $\gamma<1$ is immediate.
\end{remark}

\begin{remark}

Instead of assuming $\{(x_t,\epsilon_t)\}$ is jointly Subweibull, it is sufficient to assume the marginals $\{x_t\}$ and $\{\epsilon_t\}$ are Subweibull (see \citep{wong2020lasso,kuchibhotla2018moving}). Let the Subweibull norms of $\{x_t\}$ and $\{\epsilon_t\}$ be bounded above by $K_x$, $K_\epsilon$ respectively. Let $v\in\mbB_2$ be split as $v=(v_1,v_2)$ with $v_1\in\mbR^{p-1}$. Then

\begin{align}
||v^\top(x,\epsilon)||_{\psi_{\gamma_2}} 
& =||v_1^\top x+v_2\epsilon||_{\psi_{\gamma_2}} \\
& \overset{(a)}{\le} 2^{1/\gamma_2} \parenf{||v_1^\top x||_{\psi_{\gamma_2}}+ ||v_2\epsilon||_{\psi_{\gamma_2}} } \\
& \overset{(b)}{\le} 2^{1/\gamma_2}\parenf{||v_1|| ||x||_{\psi_{\gamma_2}}+||v_2|||\epsilon||_{\psi_{\gamma_2}} } \\
& \le 2^{1/\gamma_2} \parenf{||x||_{\psi_{\gamma_2}}+ ||\epsilon||_{\psi_{\gamma_2}} } .
\end{align}

We used (a) Lemma A.3, \citep{gotze2019concentration}, and (b) lemma 12, \citep{wong2020lasso}. Hence, $K\le2^{1/\gamma_2}(K_x +K_\epsilon)$. Also, it is equivalent to assume the Subweibull norms of $\{x_t\}$ and response $\{y_t\}$ are bounded by $K_x$ and $K_y$ (say), since, by linearity,  $||\epsilon||_{\psi_{\gamma_2}}\le ||y||_{\psi_{\gamma_2}} + ||x^\top\beta^*||_{\psi_{\gamma_2}}\le K_y+ ||\beta^*||K_x$. Since we actually observe the samples $\{(x_t,y_t)\}$, $1\le t\le n$, we might as well make assumptions on the same. Note that, in this case the linear factor is $||\beta^*||=\mcO(\sqrt{s})$, if $\beta^*$ is $s$-sparse.

\end{remark}

\begin{remark}

The pair $(\gamma_1,\gamma_2)$ encompasses a whole family of problems with the first coefficient measuring the strength of temporal dependence, and the second, measuring the heaviness of the tails of the stochastic processes. As noted in \citep{wong2020lasso}, the challenging cases are the ones where $\gamma_1 < 1$ signifying strong temporal dependence and $\gamma_2 < 2$ signifying tails heavier than a Gaussian distribution. For the  usual $\ell_1$ LASSO with independent Subweibull elements (i.e. $\gamma_1=\infty$), see Theorem 4.5, \citep{kuchibhotla2018moving}.

\end{remark}

\begin{remark}

The sub-geometric index $\gamma_1$ and the the Subweibull tail index $\gamma_2$ do not depend on underlying model parameters. However, the mixing rate $\beta_{mix}$ of the joint process $\{(x_t,\epsilon_t)\}$ often $\textit{does}$ (endogeneity). This is often an unavoidable feature of time series in general.

\end{remark}

We can estimate the regression parameter $\beta^*$ using either LASSO or the Dantzig selector as in the Gaussian case. We only need to show that the deviation and RE conditions still hold with high probability. However, unlike the Gaussian case, the temporal dependence factors that appear in the following propositions and proofs, will be a function of $\beta_{mix}$, $\gamma_1$ and $\gamma_2$, but the exact form is not explicit. This is clearly undesirable, but we feel it is unavoidable in the current setting, since we borrow our concentration results (Proposition \eqref{ConcWong}) from \citep{wong2020lasso}, which suffers from the same drawback. On the other hand, to the best of our knowledge, no explicit forms can be derived in general, without further simplifying assumptions.

\begin{proposition}{Deviation Condition for Heavy Tails.}\label{dev-Subweibull}
There is an absolute constant $c>0$, and a constant $\msfC_{dev}>0$ depending on $\beta_{mix}$, $\gamma_1$ and $\gamma_2$ such that, for $n\ge n_{dev} :=\parenf{\msfC_{dev}^2 c w^2[\mbB_\mcR(0,1)]}^{2/\gamma -1}$, we have
\begin{gather}
\mbP\left[ \mcR^*\left( \frac{u^\top X^\top\mcE}{n} \right) >\bar{\Phi}_\mcR K^2\msfC_{dev}\sqrt{ \frac{cw^2[\mbB_\mcR(0,1)]}{n}}\right] \\
\le\exp\parens{-w^2(\mbB_\mcR(0,1))+\log n} .
\end{gather}    
\end{proposition}

\begin{remark}
The results are remarkably similar to Proposition \eqref{dev-stochreg}. The only major difference is the exponent term $2/\gamma-1$ appearing in the lower bound for the sample size. In the best case scenario, when the dependence is negligible ($\gamma_1 \rightarrow\infty$) and the tails are close to being sub-Gaussian ($\gamma_2\rightarrow 2$), we have $\gamma\rightarrow 1$, so that the exponent $2/\gamma-1 \rightarrow 1$. Hence the sample size scales just as in the independent sub-Gaussian case. However, the situation gets worse when $\gamma\rightarrow 0$.  
\end{remark}

\begin{remark}
The `$\log n$' term appearing in the right hand side of the last inequality is a pretty reasonable price to pay; for example, with usual $\ell_1$ norm, i.e. $\mcR()=||\cdot||_1$, the quantity $w^2(\mbB_\mcR(0,1))$ is $\mathcal{O}(\log p)$, so with large enough constants, we have that the RHS is $\mathcal{O}(p^{-c})$ assuming $p\gg n$, which is the usual error bound for the $\ell_1$-LASSO.
\end{remark}

Next, we have to verify that the RE condition holds with high probability. So, we have the following proposition. (Again, the proof is similar to that of Proposition \eqref{RE-stochreg}. The only difference is showing Step 1 of that proof, i.e. a single concentration bound.) To fix ideas, let us consider the Dantzig selector--- the anti-concentration cone for the RE condition is $\mcT$ in this case (see Proposition \eqref{RE-stochreg}). 

\begin{proposition}{RE Condition for Heavy Tails.}\label{RE_condition_Subweibull}
Assume $\Lambda_{\min}(\Sigma_x) >0$. Then, there is a absolute constant $c'>0$ and a constant $\msfC_{RE}>0$ depending on $\beta_{mix}$, $\gamma_1$, $\gamma_2$, such that, for a minimum sample size 
\begin{gather}
n\ge n_{RE}:=\parenf{\msfC_{RE}^2 c'\max\parenc{1,\frac{16\bar{\Phi}_{\mcR}^2 K^4}{\Lambda_{\min}^2(\Sigma_X)}} w^2[\mcT\cap\mbB_2)}^{1/\gamma},    
\end{gather}
we have 
\begin{gather}
\mbP\parens{\underset{v\in\mcT\cap\mbB_2}{\inf}\frac{v^\top X^\top X v}{n}\ge\alpha_{RE} }\ge 1-\exp[-w^2(\mcT\cap\mbB_2) +\log n],
\end{gather}
where the restricted eigenvalue is $\alpha_{RE} =\Lambda_{\min}(\Sigma_x)/2$.

\end{proposition}

\begin{remark}
Note that, similar to the deviation condition \eqref{dev-Subweibull}, the only change is the exponent $1/\gamma$ appearing in the expression for the minimum sample size requirement. Again, when we approach independence ($\gamma_1\rightarrow\infty$) and sub-gaussian tails $\gamma_2\rightarrow 2$, so that the exponent $1/\gamma\rightarrow 1$, the minimum sample size for RE matches Proposition \eqref{RE-stochreg}.
\end{remark}

\begin{proposition}\label{det-errbd-SW-stochreg}{Theoretical consistency for Subweibull tails.}
Suppose the sample size and tuning parameter satisfies 
\begin{gather}
n\ge\max\parenc{{n_{dev},n_{RE}}} \quad
\lambda_n=\bar{\Phi}_\mcR K^2\msfC_{dev} \sqrt{\frac{cw^2[\mbB_\mcR(0,1)]}{n}}.  
\end{gather}
Then, the penalized estimate $\hat{\beta}$ satisfies:
\begin{gather}
||\hat{\beta}-\beta^*||\leq \frac{2\lambda_n\Phi_\mcR(\mcT)}{\alpha_{RE}}, \quad
\mcR(v)\le\frac{2\lambda_n\Phi_\mcR^2(\mcT)}{\alpha_{RE}},
\quad\text{(Estimation error)}, \\
\frac{v^\top X^\top Xv}{n}\le \frac{4\lambda_n^2\Phi_\mcR(\mcT)}{\alpha_{RE}} 
\quad\text{(Prediction error)}.
\end{gather} 
with probability at least $1-n\exp[-w^2(\mbB_\mcR(0,1))]-n\exp[-w^2(\mcT\cap\mbB_2)]$.

\end{proposition}

\begin{remark}
Thus, the message is that, in presence of temporal dependence and heavy tails, the price we pay in terms of sample size $n$ is an exponent of $\max\parenc{2/\gamma-1, 1/\gamma} =2/\gamma-1$ (since $\gamma<1$ by assumption) which, in the best case scenario is close to 1 (as $\gamma\rightarrow 1$) and becomes arbitrarily large in the worst case scenario (when $\gamma\rightarrow 0$). The order of consistency is, however, same as the Gaussian case. One issue is that the temporal dependence factors $\msfC_{dev}$ and $\msfC_{RE}$ are hard to quantify unlike the Gaussian case, where they can be written in terms of the ``spike'' $\mcM(f_X)$ of the spectral density $f_X$. However, if we have a simple model like a linear time series, for example, then the dependence factor is easy to quantify, as we see in the case of Subweibull VAR models.
\end{remark}

\section{Proofs of propositions.}

\subsection{Proof of Proposition \texorpdfstring{\eqref{DEandRE-SW}.}{}}
 
The proof is virtually identical to lemmas 5.1, 5.2 in \citep{zheng2019testing}. The only change is that, instead of the Hanson-Wright inequality for the subgaussian case \citep{rudelson2013hanson}, we use a version for the Subweibull case (\citep*[Proposition 1.1]{gotze2019concentration}). We reproduce the proof technique in \citep{zheng2019testing} in its full generality.

Let $X_t=\sum_{j\ge 0}\Psi_j\epsilon_{t-j}$ be a linear process in $\mathbb{R}^p$ where the noise components $\epsilon_{tj}$  are independent Subweibull random variables with $\sw{\epsilon_{tj}}\le \tau$. Suppose we have data $\{X_0,\cdots,X_T\}$. Let $B$ be a $p\times p$ symmetric matrix (the symmetry is required)---we consider the chaos term $\frac{1}{T}\sum_{t=0}^{T-1} X_t^\top B X_t$. In our example, we need a concentration bound for the term 
\begin{gather}
u^\top\parens{\frac{1}{T}\sum_{t=0}^{T-1} X_t X_t^\top}u
=\frac{1}{T}\sum_{t=0}^{T-1} X_t^\top uu^\top X_t,
\end{gather}
so that $B=u u^\top$ in our example. 

The chaos term $\frac{1}{T}\sum_{t=0}^{T-1} X_t^\top B X_t$ can be broken into 3 parts, as per lemma 5.2 in \citep{zheng2019testing}. We have
\begin{gather}
\begin{split}
\frac{1}{T}\sum_{t=0}^{T-1}X_t^\top B X_t=& \frac{1}{T}\sum_{t=0}^{T-1}\left(\sum_{j=0}^{\infty}\Psi_j\epsilon_{t-j}\right)^\top B\left(\sum_{j=0}^{\infty}\Psi_j\epsilon_{t-j}\right)\\
=&\frac{1}{T}\sum_{t=0}^{T-1}\left(\sum_{j=0}^{t+m}\Psi_j\epsilon_{t-j}\right)^\top B\left(\sum_{j=0}^{t+m}\Psi_j\epsilon_{t-j}\right)\\
&+\frac{1}{T}\sum_{t=0}^{T-1}\left(\sum_{j=t+m+1}^{\infty}\Psi_j\epsilon_{t-j}\right)^\top B\left(\sum_{j=t+m+1}^{\infty}\Psi_j\epsilon_{t-j}\right)\\
&+\frac{2}{T}\sum_{t=0}^{T-1}\left(\sum_{j=0}^{t+m}\Psi_j\epsilon_{t-j}\right)^\top B\left(\sum_{j=t+m+1}^{\infty}\Psi_j\epsilon_{t-j}\right)\\
\triangleq &E_1+E_2+E_3.
\end{split}    
\end{gather}
Note that, for the two cross product terms are equal because $B$ is symmetric. Here, $m$ is a positive integer to be chosen later. Then we can bound each $E_i$ from its expectation separately, and $m$ will be chosen to be sufficiently large later.

\subsubsection{Bounding \texorpdfstring{$E_1-\mathbb{E}(E_1)$:}{}}
 
Fix an index $0\le t\le T-1$. 
Let $\Theta^{(t)}\in\mathbb{R}^{p\times (T+m)p}$ and ${\epsilon}\in\mathbb{R}^{(T+m)p}$ be defined as
\begin{equation*}
\Theta^{(t)}=\begin{pmatrix}\Psi_{t+m}^{(p)}&\cdots&\Psi_0^{(p)}&0&\cdots&0\end{pmatrix}  ,
{\epsilon}=\begin{pmatrix}
\epsilon_{-m}^\top &\cdots &\epsilon_{T-1}^\top
\end{pmatrix}^\top .
\end{equation*}
Then $E_1={\epsilon}^\top \left(\frac{1}{T}\sum_{t=0}^{T-1}\Theta^{(t)\top}B\Theta^{(t)}\right){\epsilon}$, and by Proposition 1.1 in \citep{gotze2019concentration}, we only need to bound the operator norm and Frobenius norm of $\frac{1}{T}\sum_{t=0}^{T-1}\Theta^{(t)\top}B\Theta^{(t)}$.

\subsubsection{Bounding \texorpdfstring{ $\left\|\frac{1}{T}\sum_{t=0}^{T-1}\Theta^{(t)\top}B\Theta^{(t)}\right\|_2$:}{}}

For any unit vector $u,v\in \mathbb{R}^{(T+m)p}$, 
\begin{equation*}
\begin{split}
u^\top\frac{1}{T}\sum_{t=0}^{T-1}\Theta^{(t)\top}B\Theta^{(t)}v=
&\frac{1}{T}\sum_{t=0}^{T-1}\sum_{i,j=1}^{t+m+1} u^{(i)\top}\Psi_{t+m+1-i}^\top B\Psi_{t+m+1-j}v^{(j)}\\
=&\frac{1}{T}\sum_{i,j=1}^{T+m}u^{(i)\top}\left[\sum_{t=(i\vee j-m-1)\vee 0}^{T-1} \Psi_{t+m+1-i}^\top B\Psi_{t+m+1-j}\right]v^{(j)}\\
\leq & \frac{1}{T}\sum_{i,j=1}^{T+m}\|u^{(i)}\|_2\|v^{(j)}\|_2\|B\|_2\sum_{l=0}^{\infty}\left\|\Psi_{|i-j|+l}\right\|_2\left\|\Psi_l\right\|_2,
\end{split}
\end{equation*}

where $u^{(i)}=(u_{(i-1)p+1},\dots,u_{ip})$, for $1\le i\le (T+m)$ etc. Let ${\alpha}_i=\left\|\Psi_i\right\|_2$, and $\Gamma\in \mathbb{R}^{(T+m)\times (T+m)}$ be defined as
$\Gamma_{ij}=\sum_{k=0}^\infty {\alpha}_{|i-j|+k}{\alpha}_{k}$, then 
\begin{equation*}
\begin{split}
u^\top \frac{1}{T}\sum_{t=0}^{T-1}\Theta^{(t)\top} B\Theta^{(t)}v\leq 
\frac{\|B\|_2}{T}(\|u^{(1)}\|_2,\dots,\|u^{(T+m)}\|_2)\Gamma 
\begin{pmatrix}
\|v^{(1)}\|_2\\
\vdots\\
\|v^{(T+m)}\|_2
\end{pmatrix}
\leq \frac{\|B\|_2\Lambda_{\max}(\Gamma)}{T}.
\end{split}
\end{equation*}
Thus we only need to bound $\Lambda_{\max}(\Gamma)$. Applying Lemma C4 in \citep{zheng2019testing}, the largest eigenvalue of Toeplitz matrix $\Gamma$ can be bounded by
\begin{equation*}
\begin{split}
\Lambda_{\max}(\Gamma)\leq \text{ess}\sup_{\lambda}\left|\sum_{l=-\infty}^{\infty}\sum_{j=0}^{\infty}{\alpha}_{|l|+j}{\alpha}_{j}e^{il\lambda}\right|\\
\leq 2\sum_{l=0}^{\infty}\sum_{j=0}^{\infty} \alpha_{l+j}{\alpha}_{j}
= 2\msfC(\Psi).
\end{split}
\end{equation*}

So we get, 
$\left\|\frac{1}{T}\sum_{t=0}^{T-1}\Theta^{(t)\top}B\Theta^{(t)}\right\|_2 \leq \frac{2\msfC(\Psi)\|B\|_2}{T}$.
        
\subsubsection{Bounding \texorpdfstring{$\left\|\frac{1}{T}\sum_{t=0}^{T-1}\Theta^{(t)\top}B\Theta^{(t)}\right\|_{F}^2$ :}{} }

We have
\begin{equation*}
\begin{split}
\left\|\frac{1}{T}\sum_{t=0}^{T-1}\Theta^{(t)\top}B\Theta^{(t)}\right\|_{F}^2\leq \frac{1}{T^2}\sum_{s,t=0}^{T-1}\left|\text{tr}\left(\Theta^{(s)\top}B\Theta^{(s)}\Theta^{(t)\top}B\Theta^{(t)}\right)\right|,
\end{split}
\end{equation*}
Since $B$ is symmetric using its spectral decomposition $B=P^\top \Lambda P$ with orthogonal $P$ and diagonal $\Lambda$, we get
\begin{equation*}
\begin{split}
\left|\text{tr}\left(\Theta^{(s)\top}B\Theta^{(s)}\Theta^{(t)\top}B\Theta^{(t)}\right)\right|\\
\left|\text{tr}\left(P\Theta^{(s)}\Theta^{(t)\top}B\Theta^{(t)}\Theta^{(s)\top}P^\top \Lambda\right)\right|\\
\leq\|B\|_{tr}\left\|\Theta^{(s)}\Theta^{(t)\top}B\Theta^{(t)}\Theta^{(s)\top}\right\|_2\\
\leq\|B\|_{tr}\|B\|_2\left\|\Theta^{(s)}\Theta^{(t)\top} \right\|_2^2.
\end{split}
\end{equation*}

Also,
\begin{equation*}
\begin{split}
\sum_{s,t=0}^{T-1}\left\|\Theta^{(s)}\Theta^{(t)\top}\right\|_2^2 \\
=\sum_{s,t=0}^{T-1}\left\|\sum_{i=1}^{t\wedge s +m}\Psi_{t+m-i}\Psi_{s+m-i}\right\|_2^2\\
\leq \sum_{s,t=0}^{T-1}\left(\sum_{i=1}^{t\wedge s +m}{\alpha}_{t+m-i}{\alpha}_{s+m-i}\right)^2\\
=\sum_{s,t=0}^{T-1} \left(\sum_{i=0}^{t\wedge s
+m-1}\alpha_i \alpha_{|t-s|+i}\right)^2 \\
\leq \sum_{l=0}^{T-1}2(T-l)\left(\sum_{i=0}^{\infty}
\alpha_i\alpha_{l+i}\right)^2 \\
\le 2T \left[\sum_{l=0}^{T-1} \left( \sum_{i=0}^\infty \alpha_i\alpha_{l+i} \right) \right]^2
\le 2T\msfC^2(\Psi)
\end{split}
\end{equation*}

Hence we get $\left\|\frac{1}{T}\sum_{t=0}^{T-1}
\Theta^{(t)\top} B\Theta^{(t)}\right\|_{F}^2\leq \frac{2\msfC^2(\Psi)\|B\|_2\|B\|_{tr}}{T}$. By Proposition 1.1 in \citep{gotze2019concentration}, we arrive at
\begin{equation*}
\mathbb{P}\left(\left|E_1-\mathbb{E}(E_1)\right|>\delta\right)\leq 2\exp\left[-c\min\left\{ \left(\frac{T\delta}{\tau^2\msfC(\Psi)\|B\|_2}\right)^{\alpha/2}, \frac{T\delta^2}{\tau^4\msfC^2(\Psi)\|B\|_2\|B\|_{tr}}
\right\}\right].
\end{equation*}
    
\subsubsection{Bounding \texorpdfstring{ $E_2-\mathbb{E}(E_2)$:}{}}

We will show that $\left|E_2-\mathbb{E}(E_2)\right|$ is a Subweibull random variable whose Subweibull norm can be bounded above when $m$ is large enough. First we bound $\sw[\alpha/2]{E_2}$. We have
\begin{equation*}
\begin{split}
\sw[\alpha/2]{E_2} \leq \frac{2^{1+\alpha/2}\|B\|_2}{T}\sum_{t=0}^{T-1}\left(\sum_{j=t+m+1}^\infty \alpha_j\sw{||\epsilon_{t-j}||_2}\right)^2\\
\leq \frac{2^{1+\alpha/2}||B||_2
\sw{||\epsilon_0||_2}}{T}\sum_{t=0}^{T-1}
\left(\sum_{j=t+m+1}^\infty \alpha_j\right)^2\\
\leq 2^{1+\alpha/2}||B||_2
\sw{||\epsilon_0||_2} \left(\sum_{j=m}^\infty \alpha_j\right)^2.
\end{split}
\end{equation*}
By corollary A.5 in \citep{gotze2019concentration}, we have, for some constant $c(\alpha)$ that depends on $\alpha$, that
\begin{equation*}
\sw[\alpha/2]{E_2 -\mathbb{E}(E_2)} \le c(\alpha) \sw[\alpha/2]{E_2} \\
\le c(\alpha) 2^{1+\alpha/2}||B||_2
\sw{||\epsilon_0||_2} \left(\sum_{j=m}^\infty \alpha_j\right)^2
\le \frac{\tau^2 \msfC(\Psi)||B||_2}{T}
\end{equation*}
for $m$ large enough, since stability ensures $\sum_{j\ge 0} \alpha_j <\infty$ (hence the tail of this series converges to 0). Note that $\sw{||\epsilon_0||_2}$ may grow with the dimension $p$, however, in the non-asymptotic framework, the dimension $p$ is large but finite, and $m$ is completely in our control, so we will choose $m$ to be correspondingly large enough. Finally by definition of a Subweibull random variable, we have    

\begin{equation*}
\mathbb{P}\left(\left|E_2-\mathbb{E}(E_2)\right|>\delta\right)\leq 2\exp\left[-c\left(\frac{\delta T}{\tau^2\msfC(\Psi)\|B\|_2}\right)^{\alpha/2}\right].
\end{equation*}
 
\subsubsection{Bounding \texorpdfstring{ $E_3-\mathbb{E}(E_3)$:}{}}

The term $E_3$ is similar to $E_2$ in that it is also a "cross-product" of two sums, the difference being that one of the sums is finite. Hence, one can bound $E_3$ similarly by noting that $\mathbb{E}(E_3)=0$ and adapting Lemma 6 in \citep{wong2020lasso}, that $\sw[\alpha/2]{XY}\le 2^{2/\alpha} \sw{X}\sw{Y}$ for two Subweibull random variables $X$ and $Y$ (in order to deal with the cross product as before). 

In conclusion, for any $\delta>0$,
\begin{equation*}
\begin{split}
\mathbb{P}\left(\left|\frac{1}{T}\sum_{t=0}^{T-1} X_t^\top B X_t-\text{tr}(B\Sigma_X(0))\right| >\delta\right) 
\leq\sum_{i=1}^3 \mathbb{P}\left(|E_i-\mathbb{E}(E_i)|> \delta/3\right) \\
\leq 6\exp\left[-c\min\left\{\left(
\frac{T\delta}{\tau^2 \msfC(\Psi) ||B||_2}\right)^{\alpha/2},
\frac{T\delta^2}{\tau^4\msfC^2(\Psi)\|B\|_2\|B\|_{tr}}
\right\} \right].
\end{split}
\end{equation*}

The proof of the deviation condition, is similar, and mimics \citep*[lemma 5.1]{zheng2019testing}.

We point out that there exists a ``sketch'' of proving concentration inequalities for linear Subweibull processes in \citep*[Appendix C]{lin2020regularized}. However, the authors appear to make a serious mistake in their reasoning: they assume that they can rotate a random vector to make its components independent: if a random vector $x$ in $\mbR^p$ satisfies $\mbE(x)=0$, $\mbE(xx^\top)=\Sigma$, then $y=\Sigma^{-\frac{1}{2}}x$ satisfies  $\mbE(y)=0$, $\mbE(yy^\top)=I_p$. However, $\textit{this does not imply that the components of y are independent}$. This works most notably when $x$ is a Gaussian vector (this is just the trick used for proving \citep*[Proposition 2.4]{basu2015regularized}, which the authors adapt). However, if $x$ is $\textit{not}$ Gaussian, this trick fails, even for the subgaussian case (this was pointed out precisely in \citep{zheng2019testing}). Hence, a standard Hanson-Wright inequality for independent random variables cannot be applied directly.

\subsection{Proof of Proposition \texorpdfstring{\eqref{dev-stochreg}}{}.}

\begin{proof}

We divide the proof into the following steps:

Step 1: Variational characterization of the norm : Note that 
\begin{gather}
\mcR^*\parenf{\frac{X^\top\mcE}{n}} = \underset{u \in \mbB_\mcR(0,1) }{\sup} \frac{u^\top X^\top\mcE}{n} .
\end{gather}
Step 2: Single deviation bound: Fix $u \in \mbB_\mcR(0,1)$. Then $\{u^\top x_t\}$ is a centered stationary Gaussian process with 
$f_{u^\top X} = u^\top f_X u$, and $f_{u^\top X, \epsilon} = u^\top f_{X,\epsilon}$, and, since we fixed $u \in \mbB_\mcR(0,1)$, we have $\mcM(f_{u^\top X}) \leq \bar{\Phi}_\mcR^2 \mcM(f_X)$, and $\mcM(f_{u^\top X,\epsilon} ) \leq \bar{\Phi}_\mcR \mcM(f_{X,\epsilon})$. By Proposition 2.4(b) in \citep{basu2015regularized}, we get  
\begin{gather}
\mbP[u^\top X^\top\mcE/n >2\pi(\mcM(f_{u^\top X}) + \mcM(f_\epsilon) + \mcM(f_{u^\top X, \epsilon})t)] \leq 6\exp[-c_1 n\min(t,t^2)] .
\end{gather}
Step 3: Discretize : First, $\mbB_\mcR(0,1)\subseteq \mbB_2(0,\bar{\Phi}_\mcR)=\bar{\Phi}_\mcR\mbB_2$, so $\mbB_\mcR(0,1)$ is totally bounded with respect to the usual Euclidean topology. In other words, it can be covered by finitely many $\textit{ Euclidean}$ balls of any radius $\epsilon>0$ of our choosing. Such a collection is called an $\epsilon$-net, and the  $\textit{smallest}$ cardinality of such a collection is called the $\epsilon$-covering number. Choose a finite 1/4-net N of $\mbB_\mcR(0,1)$ with the $\textit{smallest}$ cardinality (which corresponds to the 1/4-$\textit{covering}$ number of $\mbB_\mcR(0,1)$). We have, 
\begin{gather}
\underset{u \in \mbB_\mcR(0,1) }{\sup} \frac{u^\top X^\top\mcE}{n} \leq 2\underset{u \in N }{\sup} \frac{u^\top X^\top\mcE}{n} .
\end{gather}
Step 4: Union bound: Since N is \textit{finite} with cardinality $card(N)$ (say), we have a finite union bound: %
\begin{gather}
\mbP[ \mcR^*\parenf{X^\top\mcE/n} > 4\pi(\bar{\Phi}_\mcR^2 \mcM(f_X) + \mcM(f_\epsilon)+ \bar{\Phi}_\mcR\mcM(f_{X,\epsilon}))t] \\
\le\mbP[\underset{u\in N}{\sup}\parenf{u^\top X^\top\mcE/n}> 2\pi(\bar{\Phi}_\mcR^2 \mcM(f_X) + \mcM(f_\epsilon)+ \bar{\Phi}_\mcR\mcM(f_{X,\epsilon}))t] \\
\leq 6\exp[-c_1 n\min(t,t^2) +\log card(N)] \\
\leq 6\exp[-c_1 n\min(t,t^2) + c_2 w^2(\mbB_\mcR(0,1)] ,
\end{gather}
where the last step follows from Sudakov's minoration inequality (see e.g. Corollary 7.4.3, \citep{vershynin2018high}): $\log card(N)\leq c_2 w^2(\mbB_\mcR(0,1))$.

Finally, choose $n,t>0$ so that
\begin{gather}
\min(t,t^2)=t^2,\quad 
-c_1 n\min(t,t^2)=(c_2+1)w^2(\mbB_\mcR(0,1)).    
\end{gather}
Let $c_0=(c_2+1)/c_1$. Then $n\ge c_0w^2(\mbB_\mcR(0,1))$ and $t =\sqrt{\frac{c_0w^2(\mbB_\mcR(0,1))}{n}}$. The proof is now complete.

\end{proof}

\subsection{Proof of Proposition \texorpdfstring{\eqref{RE-stochreg}}{}.}

\begin{proof}

Step 1: Single concentration bound: Fix $v\in\mcT\cap\mbB_2$. Then we may assume without loss, that $||v||=1$, since by definition of a cone, $v\in\mcT\cap\mbB_2$ if and only if, $v/||v||\in\mcT\cap\mbB_2$. By Proposition 2.4(a) in \citep{basu2015regularized} 
\begin{gather}
\mbP\left[\left|v^\top(X^\top X/n -\Sigma_X(0))v\right| \geq 2\pi\mcM(f_X)t\right]\leq 2\exp[- c_1'n\min(t,t^2)] . \\
\end{gather}
Fix some $t_0>0$ to be chosen later and note that $\min(t_0,t_0^2)\ge\min(1,t_0^2)$.

Step 2 :Discretize the spherical cap $\mcT\cap\mbB_2$. Note that $\mcT\cap\mbB_2\subseteq\mbB_2$, so it is totally bounded with respect to Euclidean topology. Using an $1/4$-net covering with the smallest cardinality, then taking a union bound, we get, as before
\begin{gather}
\mbP\left[\underset{v\in\mcT\cap\mbB_2}{\sup} \left|v^\top(X^\top X/n -\Sigma_X(0) )v\right|\geq 4\pi\mcM(f_X)t_0\right] \\
\leq 2\exp[- c_1'n\min(1,t_0^2)+ c_2'w^2(\mcT\cap\mbB_2)] 
\leq 2\exp[-w^2(\mcT\cap\mbB_2)].
\end{gather}
The last inequality holds whenever $c_1'n\min(1,t_0^2)\ge (c_2'+1)w^2(\mcT\cap\mbB_2)$, that is, the sample size $n$ must satisfy $n\ge c_0'\max(1,t_0^{-2}) w^2(\mcT\cap\mbB_2)$, where $c_0'=(c_2'+1)/c_1'$.

Step 3: Decenter the quadratic form: By Proposition 2.3 in \citep{basu2015regularized}, we have $v^\top\Sigma_X(0)v\ge \Lambda_{\min}(\Sigma_X(0))\ge 2\pi\mfm(f_X)$, for $||v||=1$. Thus,
\begin{gather}
\frac{v^\top X^\top X v}{n}\ge 2\pi\mfm(f_X)-4\pi\mcM(f_X)t_0,
\quad \forall v\in \mcT\cap\mbB_2,    
\end{gather}
with probability at least $1-2\exp[-w^2(\mcT\cap\mbB_2)]$. Finally, we set $t_0=\frac{\mfm(f_X)}{4\mcM(f_X)}<1$ and complete the proof.

\end{proof}

\subsection{Proof of Proposition \texorpdfstring{\eqref{det-errbd-dant-stochreg}}{}.}

\begin{proof}

Step 1: Since $\beta^*$ is feasible, we have
\begin{gather}
\mcR^*\parenf{\frac{X^\top(y - X \beta^*)}{n}} \leq \lambda_n , \\
\text{i.e.} \quad \mcR^*\parenf{\frac{X^\top\mcE}{n}}\le\lambda_n.
\end{gather} 
Step 2: Basic Inequality: Since $\hat{\beta}$ is optimal and $\beta^*$ is feasible,
\begin{gather}
\mcR(\hat{\beta}) \leq \mcR(\beta^*) \\
\Rightarrow \mcR(\beta^* + v) \leq \mcR(\beta^*) \Rightarrow v \in \mcT_\mcR(\beta^*) 
\end{gather} 
where $\mcT_\mcR(\beta^*)$ is the tangent cone at $\beta^*$.

Step 3: RE condition : Now, $v\in\mcT$,if and only if $\frac{v}{||v||}\in\mcT\cap\mbB_2$. Hence,
\begin{gather}
\alpha_{RE}||v||^2\le\frac{1}{n}||Xv||^2 = \frac{v^\top X^\top Xv}{n} \le \mcR(v) \mcR^*\parenf{\frac{X^\top Xv}{n}} . \quad \text{(H\"{o}lder's Inequality)}
\end{gather}
Step 4: By triangle inequality, we have 
\begin{gather}
\mcR^*\parenf{\frac{X^\top Xv}{n}} = \mcR^*\parenf{\frac{X^\top X(\hat{\beta} - \beta^*)}{n}} \\
\le \mcR^*\parenf{\frac{X^\top (y-X\hat{\beta})}{n}} + \mcR^*\parenf{\frac{X^\top (y-X\beta^*))}{n}}\le 2\lambda_n,
\end{gather}
which finally implies that
\begin{gather}
\alpha_{RE}||v||^2\le \frac{v^\top X^\top Xv}{n}\le 2\lambda_n  \mcR(v) \Rightarrow ||v||\le\frac{2\lambda_n \Phi_\mcR(\mcT)} {\alpha_{RE}} .
\end{gather}
It also follows 
\begin{gather}
\mcR(v)\le\Phi_\mcR(\mcT)||v||\le \frac{2\lambda_n \Phi_\mcR^2(\mcT)}{\alpha_{RE}},     
\end{gather}
which in turn implies
\begin{gather}
\frac{v^\top X^\top Xv}{n}\le 2\lambda_n\mcR(v) \le\frac{4\lambda_n^2\Phi_\mcR(\mcT)}{\alpha_{RE}} .
\end{gather}
Using the probability bounds from the deviation and RE conditions now completes the proof.

\end{proof}

\subsection{Proof of Proposition \texorpdfstring{\eqref{det-errbd-LASSO-stochreg}}{}.}

\begin{proof}

 We break up the proof in the following easy steps:

Step 1: Basic Inequality: We use the definition of $\hat{\beta}$ and (\eqref{stoch-reg-matrix-form}) to get
\begin{gather}
 \frac{1}{n} ||y- X\hat{\beta}||^{2} + \lambda_n \mcR(\hat{\beta} ) \leq  \frac{1}{n} || y- X\beta^*||^{2} + \lambda_n \mcR(\beta^*) \\
 \Rightarrow  \frac{1}{n} ||Xv||^{2} \leq \frac{2}{n} v^\top X^\top\mcE + \lambda_n[\mcR(\beta^*)- \mcR(\beta^*+v)]  
 \quad ( \text{where}\quad v=\hat{\beta} -\beta^*) \\
 \Rightarrow \frac{1}{n} ||Xv||^2 \leq \mcR(v)\mcR^*\parenf{\frac{X^\top\mcE}{n}} + \lambda_n[\mcR(\beta^*) - \mcR(\beta^*+v) ] . \quad \text{(H\"{o}lder's Inequality)} \\
 \end{gather}
We note that a general version of H\"{o}lder's inequality follows from the definition of the dual norm $\mcR^*$ : for $u,v\in\mbR^p$ we have $|u^\top v|\le \mcR(u)\mcR^*(v)$.

Step 2 : Deviation bound: Using the deviation condition stated in the proposition, we get, using L.H.S and R.H.S to denote the left and right hand sides of the last inequality:
\begin{gather}
R.H.S \leq \lambda_n [\mcR(v)/2 +\mcR(\beta^*) - \mcR(\beta^*+v)] .
\end{gather}	
Step 3: Since the $L.H.S \ge 0$, the last step implies $v \in\msC$, where the cone $\msC$ is described in the proposition.

Step 4: RE condition: From step 3, we see that that
\begin{gather}
\alpha_{RE}||v||^2\le\frsum ||Xv||^2 \le\lambda_n\parens{\frac{\mcR(v)}{2}+\mcR(\beta^*)- \mcR(\beta^*+v)} \\
\Rightarrow\alpha_{RE}||v||^2 \le\frac{1}{n} ||Xv||^2 \le \lambda_n\frac{3\mcR(v)}{2} \quad \text{(Triangle Inequality)}  .
\end{gather}
The first inequality gives $||v||\le
\frac{3\lambda_n\Phi_\mcR(\msC)}{2\alpha_{RE}}$, where $\Phi_\mcR(\msC)$ is the subspace compatibility constant. The second inequality therefore gives $||Xv||^2/n \le \frac{9\lambda_n^2 \Phi_\mcR^2(\msC)}{2\alpha_{RE}}$. This completes the proof. 
 
\end{proof}

\subsection{Proof of Proposition \texorpdfstring{\eqref{dev-Subweibull}}{}.}

\begin{proof}

Step 1: Single concentration bound: Fix $u\in\mbB_\mcR(0,1)$. Then $\{u^\top x_t\}$ is a centered, one dimensional, strictly stationary, Subweibull($\gamma_2$) process, with $||u^\top x||_{\psi_{\gamma_2}}\le \bar{\Phi}_\mcR K$. Also,
$||\epsilon||_{\psi_{\gamma_2}}\le K$. Then 
\begin{gather}
||(u^\top x)\epsilon||_{\psi_{\gamma_2/2}}\le 2^{2/\gamma_2}||u^\top x||_{\psi_{\gamma_2}}
||\epsilon||_{\psi_{\gamma_2}}
\le 2^{2/\gamma_2}\bar{\Phi}_\mcR K^2.     
\end{gather}
Again, $\mbE[u^\top x_t\epsilon_t]=\mbE[u^\top x_t\mbE[\epsilon_t|x_t]]=0$. In other words, the cross product sequence $\{u^\top x_t\epsilon_t\}$ is a centered, strictly stationary, sub-geometrically $\beta$-mixing process with finite Subweibull($\gamma_2/2$) norm. Further, by assumption, $1/\gamma=1/\gamma_1 +1/(\gamma_2/2)= 1/\gamma_1 +2/\gamma_2$ with $\gamma<1$. Hence, for $n\ge 4$, $t\ge 1/n$, we use Proposition \eqref{ConcWong} to get:
\begin{gather}
\mbP\left[ \left|\frac{u^\top X^\top\mcE}{n}\right|>t \right] \\  
\le n\exp\parens{-\frac{1}{\msfC_{dev}^2} \min\parenc{ \parenf{\frac{nt}{K_0}}^\gamma, \parenf{\frac{nt^2}{K_0^2}} } } .   
\end{gather}
Here $K_0=\bar{\Phi}_\mcR K^2$ and $\msfC_{dev}$ depends on $\beta_{mix}$, $\gamma_1$ and $\gamma_2$.

Step 2: Discretizing $\mbB_\mcR(0,1)$ and taking union bounds as before, we get, 
\begin{gather}
\mbP\left[\mcR^*\parenf{\frac{u^\top X^\top\mcE}{n}}>t \right] \\  
\le n\exp\parens{-\frac{1}{\msfC_{dev}^2} \min\parenc{ \parenf{\frac{nt}{K_0}}^\gamma, \parenf{\frac{nt^2}{K_0^2}} } + c_0 w^2[\mbB_\mcR(0,1)] } .   
\end{gather}
Note that $c_0$ is an absolute constant.

Step 3: Choose $n$, $t$ so that 
\begin{gather}
\parenf{\frac{nt}{K_0}}^\gamma \ge
\parenf{\frac{nt^2}{K_0^2}}
=\msfC_{dev}^2(c_0+1) w^2[\mbB_\mcR(0,1)] .
\end{gather}
Let $c_1=c_0+1$. We get 
\begin{gather}
n\ge\parenf{\msfC_{dev}^2 c_1 w^2[\mbB_\mcR(0,1)]}^{2/\gamma -1} ,\quad
t=K_0\msfC_{dev}\sqrt{\frac{c_1 w^2[\mbB_\mcR(0,1)]}{n}} .
\end{gather}
Note that $n$ and $t$ must also satisfy the precondition $t>1/n$. This gives
\begin{gather}
n\ge \frac{1}{K_0^2 \msfC_{dev}^2 c_1 w^2[\mbB_\mcR(0,1)]}.    
\end{gather}
Note that the last inequality is extremely mild and satisfied easily as $w^2[\mbB_\mcR(0,1)]$ grows with dimension $p$ as $p\rightarrow\infty$, specifically as long as $w^2[\mbB_\mcR(0,1)] \succsim 1$. For example, when $\mcR$ is the $\ell_1$ norm, $w^2[\mbB_\mcR(0,1)]$ is of exact order $\log p$. Hence the precondition is satisfied if $\log p\succsim 1$.

\end{proof}

\subsection{Proof of Proposition \texorpdfstring{\eqref{RE_condition_Subweibull}}{}.}

\begin{proof}

Step 1: Single concentration bound: Fix $v\in\mcT\cap\mbB_2$. Then, without loss, $||v||=1$. We have
\begin{gather}
||v^\top x||_{\psi_{\gamma_2/2}}\le 2^{2/\gamma_2}||v^\top x||_{\psi_{\gamma_2}}^2
\le 2^{2/\gamma_2}\bar{\Phi}_\mcR K^2. \\  
||v^\top x -\mbE(v^\top x)||_{\psi_{\gamma_2/2}}\le 2||v^\top x||_{\psi_{\gamma_2}}^2\le 2^{2/\gamma_2+1}\bar{\Phi}_\mcR K^2.
\end{gather}
Again, denoting $K_0=\bar{\Phi}_\mcR K^2$, we use Proposition \eqref{ConcWong} and get:
\begin{gather}
\mbP\left[\left|v^\top\left(\frac{X^\top X}{n} -\Sigma_x \right)v \right|\geq t\right]\le
 n\exp\parens{-\frac{1}{\msfC_{RE}^2} \min\parenc{ \parenf{\frac{nt}{K_0}}^\gamma, \parenf{\frac{nt^2}{K_0^2}} } } . 
\end{gather}
The constant $\msfC_{RE}>0$ depends on $\beta_{mix}$, $\gamma_1$ and $\gamma_2$.

Step 2: Discretizing the spherical cap  $\mcT\cap\mbB_2$ and taking union bounds, we get
\begin{gather}
\mbP\left[\underset{v\in\mcT\cap\mbB_2}{\sup} \left| v^\top\parenf{\frac{X^\top X}{n} -\Sigma_x}v \right|\ge 2K_0 t\right]\le
n\exp\parens{-\frac{1}{\msfC_{RE}^2} \min\parenc{(nt)^\gamma,nt^2}
+c_0w^2(\mcT\cap\mbB_2)} \\
\le n\exp\parens{-\frac{1}{\msfC_{RE}^2} n^\gamma\min(1,t^2)
+c_0w^2(\mcT\cap\mbB_2)} . 
\end{gather}
This follows because $\min((nt)^\gamma,nt^2)\ge \min(n^\gamma,n)\min(t^\gamma,t^2)\ge n^\gamma\min(1,t^2)$. We used $n\ge 4$, $\gamma<1$, $t>0$.

Step 3: Decentering the quadratic form: Choose t and n as the following:
\begin{gather}
t=\Lambda_{\min}(\Sigma_X)/4K_0^2,
\quad \frac{1}{\msfC_{RE}^2} n^\gamma\min(1,t^2)
\ge(c_0+1)w^2(\mcT\cap\mbB_2).
\end{gather}
Let $c_1=c_0+1$. Then the minimum sample size and the restricted eigenvalue are given by 
\begin{gather}
n\ge\parenf{\msfC_{RE}^2 c_1\max\parenc{1,\frac{16K_0^2}{\Lambda_{\min}(\Sigma_X)}} w^2(\mcT\cap\mbB_2)}^{1/\gamma}, \quad
\alpha_{RE}=\Lambda_{\min}(\Sigma_X)/2. 
\end{gather}
\end{proof}

The precondition $t>1/n$ boils down to $n\ge 2/\Lambda_{\min}(\Sigma_X)$.

\subsection{Proof of Proposition \texorpdfstring{\eqref{det-errbd-SW-stochreg}}{}.}

The proof is virtually identical to that of Proposition \eqref{det-errbd-dant-stochreg}. 

\subsection{Proof of Proposition \texorpdfstring{\eqref{DeviationSW-VAR}}{}.}

The proof is virtually identical to that of Proposition \eqref{dev-Subweibull}, the only difference being that we apply the single concentration bound for the deviation term $\fr{X^\top\eta}$ in Proposition \eqref{DEandRE-SW}, instead of Proposition \eqref{ConcWong}.

\subsection{Proof of Proposition \texorpdfstring{\eqref{RE-SW-VAR}}{}.}

The proof is virtually identical to that of Proposition \eqref{RE_condition_Subweibull}, the only difference being that we apply the single concentration bound for the gram matrix $\fr{X^\top X}$ in Proposition \eqref{DEandRE-SW}, instead of Proposition \eqref{ConcWong}.

\subsection{Proof of Proposition \texorpdfstring{\eqref{VAR-Subweibull-consistency}}{}.}

The proof is virtually identical to that of Proposition \eqref{det-errbd-dant-stochreg}.

\end{appendices}

%give the  reference style at the end instead before...
%...begin document coz documentclass=amsrt freaks out o.w.
\bibliographystyle{plainnat}    % reference styles

\bibliography{references}

\begin{thebibliography}{70}
\providecommand{\natexlab}[1]{#1}
\providecommand{\url}[1]{\texttt{#1}}
\expandafter\ifx\csname urlstyle\endcsname\relax
  \providecommand{\doi}[1]{doi: #1}\else
  \providecommand{\doi}{doi: \begingroup \urlstyle{rm}\Url}\fi

\bibitem[Agarwal et~al.(2012)Agarwal, Negahban, and
  Wainwright]{agarwal2012noisy}
Alekh Agarwal, Sahand Negahban, and Martin~J Wainwright.
\newblock Noisy matrix decomposition via convex relaxation: Optimal rates in
  high dimensions.
\newblock \emph{The Annals of Statistics}, pages 1171--1197, 2012.

\bibitem[Argyriou et~al.(2012)Argyriou, Foygel, and Srebro]{argyriou2012sparse}
Andreas Argyriou, Rina Foygel, and Nathan Srebro.
\newblock Sparse prediction with the $ k $-support norm.
\newblock \emph{Advances in Neural Information Processing Systems},
  25:\penalty0 1457--1465, 2012.

\bibitem[Bach et~al.(2011)Bach, Jenatton, Mairal, and
  Obozinski]{bach2011optimization}
Francis Bach, Rodolphe Jenatton, Julien Mairal, and Guillaume Obozinski.
\newblock Optimization with sparsity-inducing penalties.
\newblock \emph{arXiv preprint arXiv:1108.0775}, 2011.

\bibitem[Banerjee et~al.(2014)Banerjee, Chen, Fazayeli, and
  Sivakumar]{banerjee2014estimation}
Arindam Banerjee, Sheng Chen, Farideh Fazayeli, and Vidyashankar Sivakumar.
\newblock Estimation with norm regularization.
\newblock \emph{Advances in neural information processing systems},
  27:\penalty0 1556--1564, 2014.

\bibitem[Banerjee et~al.(2015)Banerjee, Chen, Fazayeli, and
  Sivakumar]{banerjee2015estimation}
Arindam Banerjee, Sheng Chen, Farideh Fazayeli, and Vidyashankar Sivakumar.
\newblock Estimation with norm regularization.
\newblock \emph{arXiv}, pages arXiv--1505, 2015.

\bibitem[Basu and Michailidis(2015)]{basu2015regularized}
Sumanta Basu and George Michailidis.
\newblock Regularized estimation in sparse high-dimensional time series models.
\newblock \emph{The Annals of Statistics}, 43\penalty0 (4):\penalty0
  1535--1567, 2015.

\bibitem[Basu et~al.(2015)Basu, Shojaie, and Michailidis]{basu2015network}
Sumanta Basu, Ali Shojaie, and George Michailidis.
\newblock Network granger causality with inherent grouping structure.
\newblock \emph{The Journal of Machine Learning Research}, 16\penalty0
  (1):\penalty0 417--453, 2015.

\bibitem[Basu et~al.(2019{\natexlab{a}})Basu, Das, Michailidis, and
  Purnanandam]{basu2019system}
Sumanta Basu, Sreyoshi Das, George Michailidis, and Amiyatosh Purnanandam.
\newblock A system-wide approach to measure connectivity in the financial
  sector.
\newblock \emph{Available at SSRN 2816137}, 2019{\natexlab{a}}.

\bibitem[Basu et~al.(2019{\natexlab{b}})Basu, Li, and Michailidis]{basu2019low}
Sumanta Basu, Xianqi Li, and George Michailidis.
\newblock Low rank and structured modeling of high-dimensional vector
  autoregressions.
\newblock \emph{IEEE Transactions on Signal Processing}, 67\penalty0
  (5):\penalty0 1207--1222, 2019{\natexlab{b}}.

\bibitem[Bhaskar et~al.(2013)Bhaskar, Tang, and Recht]{bhaskar2013atomic}
Badri~Narayan Bhaskar, Gongguo Tang, and Benjamin Recht.
\newblock Atomic norm denoising with applications to line spectral estimation.
\newblock \emph{IEEE Transactions on Signal Processing}, 61\penalty0
  (23):\penalty0 5987--5999, 2013.

\bibitem[Bickel and Levina(2008)]{bickel2008covariance}
Peter~J Bickel and Elizaveta Levina.
\newblock Covariance regularization by thresholding.
\newblock \emph{The Annals of Statistics}, 36\penalty0 (6):\penalty0
  2577--2604, 2008.

\bibitem[Bickel et~al.(2009)Bickel, Ritov, and
  Tsybakov]{bickel2009simultaneous}
Peter~J Bickel, Ya’acov Ritov, and Alexandre~B Tsybakov.
\newblock Simultaneous analysis of lasso and dantzig selector.
\newblock \emph{The Annals of statistics}, 37\penalty0 (4):\penalty0
  1705--1732, 2009.

\bibitem[Bogdan et~al.(2013)Bogdan, Berg, Su, and
  Candes]{bogdan2013statistical}
Malgorzata Bogdan, Ewout van~den Berg, Weijie Su, and Emmanuel Candes.
\newblock Statistical estimation and testing via the sorted l1 norm.
\newblock \emph{arXiv preprint arXiv:1310.1969}, 2013.

\bibitem[Bogdan et~al.(2015)Bogdan, Van Den~Berg, Sabatti, Su, and
  Cand{\`e}s]{bogdan2015slope}
Ma{\l}gorzata Bogdan, Ewout Van Den~Berg, Chiara Sabatti, Weijie Su, and
  Emmanuel~J Cand{\`e}s.
\newblock Slope—adaptive variable selection via convex optimization.
\newblock \emph{The annals of applied statistics}, 9\penalty0 (3):\penalty0
  1103, 2015.

\bibitem[Bondell and Reich(2008)]{bondell2008simultaneous}
Howard~D Bondell and Brian~J Reich.
\newblock Simultaneous regression shrinkage, variable selection, and supervised
  clustering of predictors with oscar.
\newblock \emph{Biometrics}, 64\penalty0 (1):\penalty0 115--123, 2008.

\bibitem[Bradley(2005)]{bradley2005basic}
Richard~C Bradley.
\newblock Basic properties of strong mixing conditions. a survey and some open
  questions.
\newblock \emph{arXiv preprint math/0511078}, 2005.

\bibitem[Brooks and Tsolacos(2000)]{brooks2000forecasting}
Chris Brooks and Sotiris Tsolacos.
\newblock Forecasting models of retail rents.
\newblock \emph{Environment and Planning A}, 32\penalty0 (10):\penalty0
  1825--1839, 2000.

\bibitem[Candes and Tao(2007)]{candes2007dantzig}
Emmanuel Candes and Terence Tao.
\newblock The dantzig selector: Statistical estimation when p is much larger
  than n.
\newblock \emph{The annals of Statistics}, 35\penalty0 (6):\penalty0
  2313--2351, 2007.

\bibitem[Chandrasekaran et~al.(2011)Chandrasekaran, Sanghavi, Parrilo, and
  Willsky]{chandrasekaran2011rank}
Venkat Chandrasekaran, Sujay Sanghavi, Pablo~A Parrilo, and Alan~S Willsky.
\newblock Rank-sparsity incoherence for matrix decomposition.
\newblock \emph{SIAM Journal on Optimization}, 21\penalty0 (2):\penalty0
  572--596, 2011.

\bibitem[Chandrasekaran et~al.(2012)Chandrasekaran, Recht, Parrilo, and
  Willsky]{chandrasekaran2012convex}
Venkat Chandrasekaran, Benjamin Recht, Pablo~A Parrilo, and Alan~S Willsky.
\newblock The convex geometry of linear inverse problems.
\newblock \emph{Foundations of Computational mathematics}, 12\penalty0
  (6):\penalty0 805--849, 2012.

\bibitem[Chen and Banerjee(2015)]{chen2015structured}
Sheng Chen and Arindam Banerjee.
\newblock Structured estimation with atomic norms: General bounds and
  applications.
\newblock In \emph{Advances in Neural Information Processing Systems}, pages
  2908--2916, 2015.

\bibitem[Chen and Banerjee(2016)]{chen2016structured}
Sheng Chen and Arindam Banerjee.
\newblock Structured matrix recovery via the generalized dantzig selector.
\newblock \emph{Advances in neural information processing systems},
  29:\penalty0 3252--3260, 2016.

\bibitem[Cheng and Pourahmadi(1993)]{cheng1993mixing}
R~Cheng and M~Pourahmadi.
\newblock The mixing rate of a stationary multivariate process.
\newblock \emph{Journal of Theoretical Probability}, 6\penalty0 (3):\penalty0
  603--617, 1993.

\bibitem[Cushman and Zha(1997)]{cushman1997identifying}
David~O Cushman and Tao Zha.
\newblock Identifying monetary policy in a small open economy under flexible
  exchange rates.
\newblock \emph{Journal of Monetary economics}, 39\penalty0 (3):\penalty0
  433--448, 1997.

\bibitem[Davidson(1994)]{davidson1994stochastic}
James Davidson.
\newblock \emph{Stochastic limit theory: An introduction for econometricians}.
\newblock OUP Oxford, 1994.

\bibitem[Fan and Li(2001)]{fan2001variable}
Jianqing Fan and Runze Li.
\newblock Variable selection via nonconcave penalized likelihood and its oracle
  properties.
\newblock \emph{Journal of the American statistical Association}, 96\penalty0
  (456):\penalty0 1348--1360, 2001.

\bibitem[Figueiredo and Nowak(2016)]{figueiredo2016ordered}
Mario Figueiredo and Robert Nowak.
\newblock Ordered weighted l1 regularized regression with strongly correlated
  covariates: Theoretical aspects.
\newblock In \emph{Artificial Intelligence and Statistics}, pages 930--938,
  2016.

\bibitem[Figueiredo and Nowak(2014)]{figueiredo2014sparse}
Mario~AT Figueiredo and Robert~D Nowak.
\newblock Sparse estimation with strongly correlated variables using ordered
  weighted l1 regularization.
\newblock \emph{arXiv preprint arXiv:1409.4005}, 2014.

\bibitem[Geraci and Gnabo(2018)]{geraci2018measuring}
Marco~Valerio Geraci and Jean-Yves Gnabo.
\newblock Measuring interconnectedness between financial institutions with
  bayesian time-varying vector autoregressions.
\newblock \emph{Journal of Financial and Quantitative Analysis}, 53\penalty0
  (3):\penalty0 1371--1390, 2018.

\bibitem[G{\"o}tze et~al.(2019)G{\"o}tze, Sambale, and
  Sinulis]{gotze2019concentration}
Friedrich G{\"o}tze, Holger Sambale, and Arthur Sinulis.
\newblock Concentration inequalities for polynomials in
  $\backslash\alpha$-sub-exponential random variables.
\newblock \emph{arXiv preprint arXiv:1903.05964}, 2019.

\bibitem[Hamilton(2020)]{hamilton2020time}
James~Douglas Hamilton.
\newblock \emph{Time series analysis}.
\newblock Princeton university press, 2020.

\bibitem[Javanmard and Montanari(2014)]{javanmard2014confidence}
Adel Javanmard and Andrea Montanari.
\newblock Confidence intervals and hypothesis testing for high-dimensional
  regression.
\newblock \emph{The Journal of Machine Learning Research}, 15\penalty0
  (1):\penalty0 2869--2909, 2014.

\bibitem[Kock and Callot(2015)]{kock2015oracle}
Anders~Bredahl Kock and Laurent Callot.
\newblock Oracle inequalities for high dimensional vector autoregressions.
\newblock \emph{Journal of Econometrics}, 186\penalty0 (2):\penalty0 325--344,
  2015.

\bibitem[Kolmogorov and Rozanov(1960)]{kolmogorov1960strong}
Andrei~Nikolaevich Kolmogorov and Yu~A Rozanov.
\newblock On strong mixing conditions for stationary gaussian processes.
\newblock \emph{Theory of Probability \& Its Applications}, 5\penalty0
  (2):\penalty0 204--208, 1960.

\bibitem[Kuchibhotla and Chakrabortty(2018)]{kuchibhotla2018moving}
Arun~Kumar Kuchibhotla and Abhishek Chakrabortty.
\newblock Moving beyond sub-gaussianity in high-dimensional statistics:
  Applications in covariance estimation and linear regression.
\newblock \emph{arXiv preprint arXiv:1804.02605}, 2018.

\bibitem[Lin and Michailidis(2020)]{lin2020regularized}
Jiahe Lin and George Michailidis.
\newblock Regularized estimation of high-dimensional factor-augmented vector
  autoregressive (favar) models.
\newblock \emph{Journal of machine learning research}, 21\penalty0
  (117):\penalty0 1--51, 2020.

\bibitem[Loh(2017)]{loh2017statistical}
Po-Ling Loh.
\newblock Statistical consistency and asymptotic normality for high-dimensional
  robust m-estimators.
\newblock \emph{The Annals of Statistics}, 45\penalty0 (2):\penalty0 866--896,
  2017.

\bibitem[Loh(2018)]{loh2018scale}
Po-Ling Loh.
\newblock Scale calibration for high-dimensional robust regression.
\newblock \emph{arXiv preprint arXiv:1811.02096}, 2018.

\bibitem[Loh and Wainwright(2012)]{loh2012high}
Pp-Ling Loh and Martin~J Wainwright.
\newblock High-dimensional regression with noisy and missing data: Provable
  guarantees with nonconvexity.
\newblock \emph{The Annals of Statistics}, 40\penalty0 (3):\penalty0
  1637--1664, 2012.

\bibitem[L{\"u}tkepohl(2005)]{lutkepohl2005new}
Helmut L{\"u}tkepohl.
\newblock \emph{New introduction to multiple time series analysis}.
\newblock Springer Science \& Business Media, 2005.

\bibitem[Melnyk and Banerjee(2016)]{melnyk2016estimating}
Igor Melnyk and Arindam Banerjee.
\newblock Estimating structured vector autoregressive models.
\newblock In \emph{International Conference on Machine Learning}, pages
  830--839, 2016.

\bibitem[Merlev{\`e}de et~al.(2011)Merlev{\`e}de, Peligrad, and
  Rio]{merlevede2011bernstein}
Florence Merlev{\`e}de, Magda Peligrad, and Emmanuel Rio.
\newblock A bernstein type inequality and moderate deviations for weakly
  dependent sequences.
\newblock \emph{Probability Theory and Related Fields}, 151\penalty0
  (3-4):\penalty0 435--474, 2011.

\bibitem[Michailidis and d’Alch{\'e}
  Buc(2013)]{michailidis2013autoregressive}
George Michailidis and Florence d’Alch{\'e} Buc.
\newblock Autoregressive models for gene regulatory network inference:
  Sparsity, stability and causality issues.
\newblock \emph{Mathematical biosciences}, 246\penalty0 (2):\penalty0 326--334,
  2013.

\bibitem[Negahban et~al.(2012)Negahban, Ravikumar, Wainwright, and
  Yu]{negahban2012unified}
Sahand~N Negahban, Pradeep Ravikumar, Martin~J Wainwright, and Bin Yu.
\newblock A unified framework for high-dimensional analysis of $ m $-estimators
  with decomposable regularizers.
\newblock \emph{Statistical science}, 27\penalty0 (4):\penalty0 538--557, 2012.

\bibitem[Nicholson et~al.(2017)Nicholson, Matteson, and
  Bien]{nicholson2017varx}
William~B Nicholson, David~S Matteson, and Jacob Bien.
\newblock Varx-l: Structured regularization for large vector autoregressions
  with exogenous variables.
\newblock \emph{International Journal of Forecasting}, 33\penalty0
  (3):\penalty0 627--651, 2017.

\bibitem[Nicholson et~al.(2020)Nicholson, Wilms, Bien, and
  Matteson]{nicholson2020high}
William~B Nicholson, Ines Wilms, Jacob Bien, and David~S Matteson.
\newblock High dimensional forecasting via interpretable vector
  autoregression., 2020.

\bibitem[Nijs et~al.(2007)Nijs, Srinivasan, and Pauwels]{nijs2007retail}
Vincent~R Nijs, Shuba Srinivasan, and Koen Pauwels.
\newblock Retail-price drivers and retailer profits.
\newblock \emph{Marketing Science}, 26\penalty0 (4):\penalty0 473--487, 2007.

\bibitem[Ning and Liu(2017)]{ning2017general}
Yang Ning and Han Liu.
\newblock A general theory of hypothesis tests and confidence regions for
  sparse high dimensional models.
\newblock \emph{The Annals of Statistics}, 45\penalty0 (1):\penalty0 158--195,
  2017.

\bibitem[Rao et~al.(2012)Rao, Recht, and Nowak]{rao2012universal}
Nikhil Rao, Ben Recht, and Robert Nowak.
\newblock Universal measurement bounds for structured sparse signal recovery.
\newblock In \emph{Artificial Intelligence and Statistics}, pages 942--950,
  2012.

\bibitem[Rudelson et~al.(2013)Rudelson, Vershynin, et~al.]{rudelson2013hanson}
Mark Rudelson, Roman Vershynin, et~al.
\newblock Hanson-wright inequality and sub-gaussian concentration.
\newblock \emph{Electronic Communications in Probability}, 18, 2013.

\bibitem[Seth et~al.(2015)Seth, Barrett, and Barnett]{seth2015granger}
Anil~K Seth, Adam~B Barrett, and Lionel Barnett.
\newblock Granger causality analysis in neuroscience and neuroimaging.
\newblock \emph{Journal of Neuroscience}, 35\penalty0 (8):\penalty0 3293--3297,
  2015.

\bibitem[Stock and Watson(2016)]{stock2016dynamic}
James~H Stock and Mark~W Watson.
\newblock Dynamic factor models, factor-augmented vector autoregressions, and
  structural vector autoregressions in macroeconomics.
\newblock In \emph{Handbook of macroeconomics}, volume~2, pages 415--525.
  Elsevier, 2016.

\bibitem[Stucky and van~de Geer(2018)]{stucky2018asymptotic}
Benjamin Stucky and Sara van~de Geer.
\newblock Asymptotic confidence regions for high-dimensional structured
  sparsity.
\newblock \emph{IEEE Transactions on Signal Processing}, 66\penalty0
  (8):\penalty0 2178--2190, 2018.

\bibitem[Talagrand(2006)]{talagrand2006generic}
Michel Talagrand.
\newblock \emph{The generic chaining: upper and lower bounds of stochastic
  processes}.
\newblock Springer Science \& Business Media, 2006.

\bibitem[Tibshirani(1996)]{tibshirani1996regression}
Robert Tibshirani.
\newblock Regression shrinkage and selection via the lasso.
\newblock \emph{Journal of the Royal Statistical Society: Series B
  (Methodological)}, 58\penalty0 (1):\penalty0 267--288, 1996.

\bibitem[van~de Geer(2014)]{van2014weakly}
Sara van~de Geer.
\newblock Weakly decomposable regularization penalties and structured sparsity.
\newblock \emph{Scandinavian Journal of Statistics}, 41\penalty0 (1):\penalty0
  72--86, 2014.

\bibitem[van~de Geer and Lederer(2013)]{van2013bernstein}
Sara van~de Geer and Johannes Lederer.
\newblock The bernstein--orlicz norm and deviation inequalities.
\newblock \emph{Probability theory and related fields}, 157\penalty0
  (1-2):\penalty0 225--250, 2013.

\bibitem[Van~de Geer et~al.(2011)Van~de Geer, B{\"u}hlmann, and
  Zhou]{van2011adaptive}
Sara Van~de Geer, Peter B{\"u}hlmann, and Shuheng Zhou.
\newblock The adaptive and the thresholded lasso for potentially misspecified
  models (and a lower bound for the lasso).
\newblock \emph{Electronic Journal of Statistics}, 5:\penalty0 688--749, 2011.

\bibitem[Van~de Geer et~al.(2014)Van~de Geer, B{\"u}hlmann, Ritov, and
  Dezeure]{van2014asymptotically}
Sara Van~de Geer, Peter B{\"u}hlmann, Ya’acov Ritov, and Ruben Dezeure.
\newblock On asymptotically optimal confidence regions and tests for
  high-dimensional models.
\newblock \emph{The Annals of Statistics}, 42\penalty0 (3):\penalty0
  1166--1202, 2014.

\bibitem[Vershynin(2018)]{vershynin2018high}
Roman Vershynin.
\newblock \emph{High-dimensional probability: An introduction with applications
  in data science}, volume~47.
\newblock Cambridge university press, 2018.

\bibitem[Wainwright(2019)]{wainwright2019high}
Martin~J Wainwright.
\newblock \emph{High-dimensional statistics: A non-asymptotic viewpoint},
  volume~48.
\newblock Cambridge University Press, 2019.

\bibitem[Wong et~al.(2016)Wong, Li, and Tewari]{wong2016lasso}
Kam~Chung Wong, Zifan Li, and Ambuj Tewari.
\newblock Lasso guarantees for time series estimation under subgaussian tails
  and $\beta$-mixing.
\newblock \emph{arXiv preprint arXiv:1602.04265}, 2016.

\bibitem[Wong et~al.(2020)Wong, Li, and Tewari]{wong2020lasso}
Kam~Chung Wong, Zifan Li, and Ambuj Tewari.
\newblock Lasso guarantees for $\beta $-mixing heavy-tailed time series.
\newblock \emph{Annals of Statistics}, 48\penalty0 (2):\penalty0 1124--1142,
  2020.

\bibitem[Wood(2009)]{wood2009presidential}
B~Dan Wood.
\newblock Presidential saber rattling and the economy.
\newblock \emph{American Journal of Political Science}, 53\penalty0
  (3):\penalty0 695--709, 2009.

\bibitem[Wu et~al.(2016)Wu, Wu, et~al.]{wu2016performance}
Wei-Biao Wu, Ying~Nian Wu, et~al.
\newblock Performance bounds for parameter estimates of high-dimensional linear
  models with correlated errors.
\newblock \emph{Electronic Journal of Statistics}, 10\penalty0 (1):\penalty0
  352--379, 2016.

\bibitem[Zeng and Figueiredo(2014)]{zeng2014ordered}
Xiangrong Zeng and M{\'a}rio~AT Figueiredo.
\newblock The ordered weighted $\ell_1$ norm: Atomic formulation, projections,
  and algorithms.
\newblock \emph{arXiv preprint arXiv:1409.4271}, 2014.

\bibitem[Zhang(2010)]{zhang2010nearly}
Cun-Hui Zhang.
\newblock Nearly unbiased variable selection under minimax concave penalty.
\newblock \emph{The Annals of statistics}, 38\penalty0 (2):\penalty0 894--942,
  2010.

\bibitem[Zhang and Zhang(2014)]{zhang2014confidence}
Cun-Hui Zhang and Stephanie~S Zhang.
\newblock Confidence intervals for low dimensional parameters in high
  dimensional linear models.
\newblock \emph{Journal of the Royal Statistical Society: Series B (Statistical
  Methodology)}, 76\penalty0 (1):\penalty0 217--242, 2014.

\bibitem[Zheng and Raskutti(2019)]{zheng2019testing}
Lili Zheng and Garvesh Raskutti.
\newblock Testing for high-dimensional network parameters in auto-regressive
  models.
\newblock \emph{Electronic Journal of Statistics}, 13\penalty0 (2):\penalty0
  4977--5043, 2019.

\bibitem[Zou and Hastie(2005)]{zou2005regularization}
Hui Zou and Trevor Hastie.
\newblock Regularization and variable selection via the elastic net.
\newblock \emph{Journal of the royal statistical society: series B (statistical
  methodology)}, 67\penalty0 (2):\penalty0 301--320, 2005.

\end{thebibliography}

\end{document}